\documentclass[12pt]{amsart}
\usepackage{amsmath,amscd,amssymb,amsthm,array}

\usepackage{mathtools}

\usepackage{hyperref}

\usepackage{amsmath,amssymb,mathrsfs,amsthm, tikz-cd,mathrsfs}

\usepackage{comments}
\newComments\DL{D}{red}
\newComments\Bj{Bj}{blue}
\newComments\AL{D}{red}

\newtheorem{mTheorem}{Theorem}[subsection]
\newtheorem{mProposition}[mTheorem]{Proposition}
\newtheorem{mLemma}[mTheorem]{Lemma}
\newtheorem{mCorollary}[mTheorem]{Corollary}

\newtheorem{mRemark}[mTheorem]{Remark}

\usepackage{DLdef1}
\newcommand{\Size}{\text{Size}}
\def\mmat #1,#2,#3,#4,{\text{\small\arraycolsep=3pt $
\begin{pmatrix}#1&#2\\#3&#4\end{pmatrix}$}}

\let\ssec\subsection

\renewcommand {\ssbegin}[2][*]
 {\refstepcounter{subsection}%
\if#1*
\addcontentsline{toc}{subsection}{\thesubsection.\hskip 1pc #2}%
\else
\addcontentsline{toc}{subsection}{\thesubsection.\hskip 1pc #2. #1}%
\fi
 \def \secno {\gdef \secno {}{\ssecfont
\thesubsection.\hskip 2ex}%
 }%
 \begin{#2}}

\renewcommand {\sssbegin}[2][*]
 {\refstepcounter{subsubsection}
\if#1*
\addcontentsline{toc}{subsubsection}{\thesubsubsection.\hskip 1pc #2}%
\else
\addcontentsline{toc}{subsubsection}{\thesubsubsection.\hskip 1pc #2. #1}
\fi
 \def \secno {\gdef \secno {}{\ssecfont \thesubsubsection.\hskip 2ex}%
 }%
 \begin{#2}}

\renewcommand {\parbegin}[2][*]
 {\refstepcounter{paragraph}
\if#1*
\addcontentsline{toc}{paragraph}{\theparagraph.\hskip 1pc #2}%
\else
\addcontentsline{toc}{paragraph}{\theparagraph.\hskip 1pc #2. #1}
\fi
 \def \secno {\gdef \secno {}{\ssecfont \theparagraph.\hskip 2ex}%
 }%
 \begin{#2}}

\begin{document}



\title[Non-degenerate anti-symmetric bilinear forms and Manin Triples] {Manin triples and non-degenerate anti-symmetric bilinear forms on Lie superalgebras in characteristic $2$}

\author{Sa\"id Benayadi}
\address{Laboratoire de Math\'ematiques IECL UMR CNRS 7502,
Universit\'e de Lorraine, 3 rue Augustin Fresnel, BP 45112,
F-57073 Metz Cedex 03, FRANCE.}

\email{said.benayadi@univ-lorraine.fr }

\author{Sofiane Bouarroudj}

 \address{Division of Science and Mathematics, New York University Abu Dhabi, Po Box 129188, Abu Dhabi, United Arab Emirates.}

\email{sofiane.bouarroudj@nyu.edu}

\thanks{The second author was supported by the grant NYUAD-065.}

\begin{abstract}
In this paper, we introduce and develop the notion of a~Manin triple for a~Lie superalgebra $\fg$ defined over a~field of characteristic $p=2$. We find cohomological necessary conditions for the pair $(\fg, \fg^*)$ to form a~Manin triple. We introduce the concept of Lie bi-superalgebras for $p=2$ and establish a~link between Manin triples and Lie bi-superalgebras.  In particular, we study Manin triples defined by a classical $r$-matrix with an extra condition (called an {\it admissible} classical $r$-matrix).   A particular case is examined where $\fg$ has an even invariant non-degenerate bilinear form.  In this case,  admissible $r$-matrices can be obtained inductively through the process of double extensions.   In addition, we introduce the notion of double extensions of Manin triples, and show how to get a new Manin triple from an existing one.
%
%
%
%
\smallskip

\end{abstract}

\maketitle

\tableofcontents
\setcounter{tocdepth}{3}

\section{Introduction} \label{SecDef}
Hereafter, $\Bbb K$ is an arbitrary field of characteristic $p$. We will specify when $\Bbb K$ is required to be algebraically closed.  For a~list of simple Lie superalgebras and their ``relatives'' we use in this paper, see \cite{BGL, BGLLS}. 
\ssec{Manin triples for Lie (super)algebras} 
Drinfeld \cite{D} proved that a~complex Lie algebra $\fg$ is also a~Lie bialgebra if and only if it constitutes, together with its dual $\fg^*$, a~Manin triple $(\fh,\fg,\fg^*)$, i.e., $\fh = \fg \oplus \fg^*$ as vector spaces, $\fh$ is a~Lie algebra having a~non-degenerate invariant symmetric bilinear form $\mathscr B$  (NIS for short) such that $\fg$ and $\fg^*$ are isotropic with respect to $\mathscr B$.  A superization of Drinfeld's result is due to Olshansky \cite{O}.  Further, Leites and Shapovalov \cite{LSh} listed, \`a la Olshlansky,  important examples of simple (and close to simple) Lie superalgebras $\fh$ of finite dimension and of polynomial growth such that $\fh$ splits into Manin-Olshansky triple.

In this paper we construct matched pairs as well as Manin triples for Lie superalgebras in characteristic $2$. Our construction is different from the case where $p\not = 2$ because of the squaring,  which was equivalent to the bracket for $p\not =2$ but plays an important role now, and introduce several completely new concepts.
\subsection{$r$-Matrices on (super)algebras}
Belavin and Drinfeld \cite{BD1, BD2} classified $r$-matrices (that describe quasi- triangular Lie bialgebra structures). Their classification uses the notion of admissible triples that satisfy certain relations.  A superization of \cite{BD1, BD2} is due to Leites and Serganova  \cite{LS} with corrections in \cite{Lyb} and by Karaali \cite{Ka1, Ka2} who used similar constructions using ideas of Belavin and Drinfeld for Lie superalgebras with a~NIS.  Karaali  constructed several examples of $r$-matrices on simple Lie superalgebras; however, the classification is not reached yet since certain non-degenerate $r$-matrices do not fit such a~description.

For Lie algebras,  Michaelis \cite{M} found a~condition under which the structure of a~triangular, coboundary Lie {\it bialgebra} exists on any Lie algebra containing linearly independent elements $a$ and $b$ satisfying $[a,b]=kb$ for some non-zero $k\in \mathbb{K}$. This method gives the structure of a~Lie bialgebra on the Virasoro algebra (with or without a~central element) and on many other Lie algebras.

Feldvoss in \cite{F1} showed that for every finite-dimensional Lie algebra over {\it an algebraically closed} field of characteristic $p\not =2$, the following are equivalent: 

(i) the algebra is either abelian or the three dimensional Heisenberg algebra; 

(ii) it does not admits any non-trivial triangular Lie bialgebra structure; 

(iii) it does not admit any non-trivial quasi-triangular Lie bialgebra structure. 

However, if the characteristic is equal to 2, then any non-abelian Lie algebra admits a~non-trivial quasi-triangular Lie bialgebra structure.

A further generalization to an arbitrary field was carried out by Feldvoss himself in \cite{F2}.

De Smedt (\cite{Sm}) in 1994 (rediscovered in \cite{F2}) showed that every finite-dimensional non-abelian Lie algebra over the real or complex numbers admits a~non-trivial coboundary Lie bialgebra structure.

For Lie superalgebras over a~field of characteristic $2$,  we define the notion of  Lie bi-superalgebra and show that this notion is equivalent to the notion of a Manin triple (see \S \ref{rmatManin}).  In the case of coboundary Lie bi-superalgebras,  the associated classical $r$-matrices must have an extra condition, see Eq. (\ref{JIsqdual})  in Theorem \ref{sqgdual}.  We call such classical $r$-matrices {\it admissible}.  In the case where the Lie superalgebra has an even NIS,  we show that the existence of an admissible $r$-matrix is equivalent to the existence of an admissible Yang-Baxter operator,  see  Prop. \ref{propUZ}.  This operator is used to study the structure of Lie bi-superalgebra associated with the admissible classical $r$-matrix. Due to the distinguish role of the squaring on the Lie superalgebra in characteristic $2$, all these notions require special approaches as compared to their analogs in other characteristics. 

\ssec{Double extensions (DEs) of Lie (super)algebras with a non-degenerate invariant symmetric bilinear form (NIS)}
Medina and Revoy in \cite{MR1} introduced the notion of \textit{double extensions} for Lie algebras with NIS in characteristic zero by means of a 1-dimensional space on a~simple Lie algebra. They showed that if $\fg$ is such an~irreducible non-simple Lie algebra (of $\dim(\fg)>1$),  then it can be described inductively in terms of another Lie algebra $\fa$ with NIS; in particular, if the centre of $\fg$ is not trivial, then $\fg$ is a~double extension of a Lie algebra $\fa$ with NIS by means of a~1-dimensional space (see also \cite{FS}).  Solvable NIS-Lie algebras $\fg$ can be embraced by this method since the centre of $\fg$ is non-trivial.  Passage to the ground field of characteristic $p>0$ brings new examples (see \cite{BBH}),  but nothing conceptually new when it comes to the construction of double extensions, except, perhaps, for inductive descriptions \`a la Medina and Revoy, since Lie's theorem and the Levi decomposition do not hold true any more if $p>0$.

The first superization of results by Medina and Revoy was initiated in \cite{BeB} provided $\mathrm{dim}(\fg_\od)=2$. They showed that every non-simple NIS-Lie superalgebra can be described by successive double extensions. The construction has been improved in a~series of papers (see \cite{ABBQ, ABB,  BBB, B,Bor}).  Superization of the definition of a~double extension introduces completely new possibilities: \textbf{the bilinear form ${\mathscr B}$, the derivation ${\mathscr D}$,  the central element $x$ may be odd} provided $p({\mathscr B})=p(x)+p(x^*)$. The first of these possibilities is related with Manin-Olshansky triples of the type described in \cite{LSh}, and loops with values in $\fp\fs\fq(n)$ and $\fh^{(1)}(0|2n+1)$, see \cite{BLS}. 

The \textbf{case of characteristic $2$} is completely different and requires new definitions and methods as shown in \cite{BeBou}. 

In the last section of this paper, the techniques of double extensions will be used to describe a class of Lie superalgebras with admissible $r$-matrices.  Such Lie superalgebras carry a NIS as well as a closed $\od$-antisymmetric ortho-orthogonal form -- see \S \ref{Bilp=2} for the terminology. An adapted process of double extension is used to describe these $r$-matrices inductively,  see \S\ref{DEeven}.  Additionally, we introduce the notion of double extension for Manin triples in \S \ref{dermanin} and establish its converse.  \\

\noindent {\bf Notation.} Let $\fg$ be a~finite-dimensional Lie superalgebra, and let $\fg^*$ be its dual also endowed with a~Lie superalgebra structure. By identifying $\fg\cong\fg^{**}$ we view the map $\ad_f(-)(x)\in \fg^{**}$, for every $f\in \fg^*$ and $x\in \fg$, as an element of $\fg$ and denote it by $x\circ \ad_f$ for short.

\section{Background}
The \textit{general linear} Lie
superalgebra $\fgl(V)$  of all linear operators in the superspace $V=V_{\bar 0}\oplus V_{\bar 1}$ over the ground field $\Kee$ is denoted by
$\fgl(\Size)$, where  $\Size:=(p_1, \dots, p_{|\Size|})$ is an ordered
collection of parities of the basis vectors of $V$ for which we take only vectors \textit{homogeneous with respect to parity}, and $|\Size|:=\dim V$. The  Lie
superalgebra of all supermatrices of size $\Size$ is also denoted by 
$\fgl(\Size)$.
Usually, for the \textit{standard} (simplest from a~certain point of view) format $\Size_{st}:=(\ev,
\dots, \ev, \od, \dots, \od)$, the notation $\fgl(\Size_{st})$ is abbreviated to $\fgl(\dim V_{\bar
0}|\dim V_{\bar 1})$. Any $X\in\fgl(\Size)$ can
be uniquely expressed as the sum of its even and odd parts; in the
standard format this is the following block expression; on non-zero summands the parity is defined:
\[
X=\mmat A,B,C,D,=\mmat A,0,0,D,+\mmat 0,B,C,0,,\quad
 p\left(\mmat A,0,0,D,\right)=\ev, \; p\left(\mmat 0,B,C,0,\right)=\od.
\]
\ssec{Symmetric and anti-symmetric bilinear forms on superspaces for $p\not =2$ see \cite[Ch.1]{LSoS}}\label{Invpn2}

Let $V$ and $W$ be two superspaces defined over an arbitrary field ${\mathbb K}$ of characteristic $p\not =2$. 

Let ${\mathscr B}\in \text{Bil}(V,W)$ be a homogenous bilinear form. The Gram matrix $B=(B_{ij})$ associated to ${\mathscr B}$ is given by the formula, 
\be\label{martBil}
B_{ij}=(-1)^{p({\mathscr B})p(v_i)}{\mathscr B}(v_{i}, v_{j})\text{~~for the basis vectors $v_{i}\in V$.}
\ee
This definition can be extended by linearity to non-homogenous forms.  Moreover,  it allows us to identify a~bilinear form $B(V, W)$ with an element of $\Hom(V, W^*)$. 
Consider the \textit{upsetting} of bilinear forms
$u\colon\Bil (V, W)\tto\Bil(W, V)$ given by the formula \be\label{susyB}
u({\mathscr B})(w, v)=(-1)^{p(v)p(w)}{\mathscr B}(v,w)\text{~~for any $v \in V$ and $w\in W$.}
\ee
In terms of the Gram matrix $B$ of ${\mathscr B}$: the form
${\mathscr B}$ is  \textit{symmetric} if  and only if 
\be\label{BilSy}
u(B)=B,\;\text{ where $u(B)=
\mmat R^{t},(-1)^{p({\mathscr B})}T^{t},(-1)^{p({\mathscr B})}S^{t},-U^{t},$ for $B=\mmat R,S,T,U,$.}
\ee
Similarly, \textit{anti-symmetry} of ${\mathscr B}$ means that $u(B)=-B$.

Since $p\not =2$, observe that if ${\mathscr B}$ is symmetric,  then $U$ is zero-diagonal, and if ${\mathscr B}$ is anti-symmetric,  then $R$ is zero-diagonal.

A non-degenerate even (resp. odd) bilinear form is called {\it orthosymplectic} (resp.  {\it periplectic}).  

Most popular normal shapes of the even non-degenerate (anti)-symmetric forms are the following matrices in the standard format:\index{$J_{2n}$}\index{$\Pi_{n}$}\index{$B_{ev}(2k\vert n)$, the normal shape of the even form}
\be\label{NorEv}
\begin{array}{l}B'_{ev}(m|2n)= \mmat 1_m,0,0,J_{2n},,\quad \text{where
$J_{2n}=\mmat 0,1_n,-1_n,0,$}\\
B_{ev}(2k|2n)= \mmat \Pi_{2k},0,0,J_{2n},, \quad \text{where
$\Pi_{2k}=\mmat 0,1_k,1_k,0,$},\\
B_{ev}(2k+1|2n)= \mmat
\Pi_{2k+1},0,0,J_{2n},, \quad \text{where
$\Pi_{2k+1}=\mat{0&0&1_k\\
0&1&0\\
1_k&0&0}$}.
\end{array}
\ee
The normal shapes of the odd non-degenerate bilinear\index{$B_{odd}(n\vert n)$, the normal shape of the odd form} form are the following matrices in the standard format (NB: we did not make a~mistake below selecting (anti-)symmetric form)
\be\label{NorOd}
B_{odd}(n|n)=\begin{cases}J_{2n}&\text{if the form is symmetric},\\
\Pi_{2n}&\text{if the form is anti-symmetric}.\\
\end{cases}
\ee 

\ssec{Symmetric and anti-symmetric bilinear forms on superspaces for $p=2$}\label{Bilp=2}
Let ${\mathscr B}\in \text{Bil}(V,W)$ be a bilinear form. The Gram matrix $B=(B_{ij})$ associated to ${\mathscr B}$ is given by the formula
\be\label{martBil}
B_{ij}= {\mathscr B}(v_{i}, v_{j})\text{~~for the basis vectors $v_{i}\in V$.}
\ee
The \textit{upsetting} of bilinear forms
$u\colon\Bil (V, W)\tto\Bil(W, V)$ is given by the formula \be\label{susyB}
u({\mathscr B})(w, v)= {\mathscr B}(v,w)\text{~~for any $v \in V$ and $w\in W$.}
\ee
In terms of the Gram matrix $B$ of ${\mathscr B}$, we have 
\[
u(B)=
\mmat R^{t}, T^{t}, S^{t}, U^{t}, \text{ for } B=\mmat R,S,T,U,.
\] 
Our main definition is the following.  We say that the bilinear form ${\mathscr B}$ is 
\begin{itemize}
\item {\it symmetric} if $u(B)=B$;
\item $\ev$-{\it antisymmetric} (called supersymmetric in  \cite{BeBou}) if $u(B)=B$ and $U$ is zero-diagonal;
\item $\od$-{\it antisymmetric} if $u(B)=B$ and $R$ is zero-diagonal.
\end{itemize}
In the particular case of a purely even space,  we say that $B$ is {\it symmetric} if $u(B)=B$ and antisymmetric if $u(B)=B$ and $R$ is zero-diagonal (recall that in this case $S=0, T=0,$ and $U=0$.)

A non-degenerate {\bf even} (resp. {\bf odd}) bilinear form is called {\it ortho-orthogonal} (resp.  {\it periplectic}).  

A. Lebedev proved (see \cite{LeD}) that over algebraically closed fields of characteristic $2$ the equivalence classes of non-degenerate symmetric bilinear forms  ${\mathscr B}$  on the $n$-dimensional space $V$ are as follows:  

(i) for $n$ odd, there is just one equivalence class, for the normal shape of ${\mathscr B}$ one can take $1_n$. 

(ii) for $n$ even, there are two equivalence classes. For the normal shape of ${\mathscr B}$ of one class one can take $1_n$, of another class one can take $\Pi_n$.

{ \bf Following the custom in differential geometry, we will denote a $\od$-antisymmetric bilinear form by $\omega$.}

Let us reformulate the condition of $\od$-{\it antisymmetry} in terms of the bilinear form $\omega$.  An even bilinear  form is $\od$-{\it antisymmetric}  if and only if:
\begin{equation}
 \label{Axiom1} \omega (x,y)=\omega(y,x),  \text{ for all $x,y\in V_{\bar i}$ where $\bar i=\ev,\od$; and $\omega(x,x)=0$ for all $x\in V_\ev$}.
\end{equation}
An odd bilinear form is $\od$-{\it antisymmetric}  if and only if:
\begin{equation}
 \label{AxiomO1b} \omega(x,y)=\omega(y,x) \text{ for all $x,y\in \fg$}; 
\end{equation}

\ssec{Lie superalgebras in characteristic $2$.} 
Following \cite{BGL, LeD,  BGLLS, BLLSq}, a~\textit{Lie superalgebra} in characteristic $2$ is a~
superspace $\fg=\fg_\ev\oplus\fg_\od$ over $\Bbb K$ such that $\fg_\ev$ is a~Lie algebra, $\fg_\od$ is a
$\fg_\ev$-module made two-sided by symmetry, and on $\fg_\od$ a~\textit{squaring}, denoted by $s_\fg:\fg_\od\rightarrow \fg_\ev$, is
given. 
The bracket on $\fg_\ev$, as well as the action of $\fg_\ev$ on $\fg_\od$, is denoted by the same symbol $[\cdot ,\cdot]_\fg$. For any $x, y\in\fg_\od$, their bracket is then defined by
\[
[x,y]_\fg:=  s_\fg(x+y)-s_\fg(x) -s_\fg(y).
\]
The bracket is extended to non-homogenous elements by bilinearity.

The Jacobi identity involving the squaring reads as follows:
\begin{equation*}\label{JIS}
[s_\fg(x),y]_\fg=[x,[x,y]_\fg]_\fg\;\text{ for any $x \in \fg_\od$ and $y\in\fg$}.
\end{equation*}
Let ${\mathscr B}_\fg$ be a bilinear form on $\fg:=(\fg, [\cdot ,\cdot]_\fg,s_\fg)$. Following \cite{BeBou}, we say that ${\mathscr B}_\fg$ is {\it invariant} if ${\mathscr B}_\fg([x,y]_\fg,z)={\mathscr B}_\fg(x,[y,z]_\fg)$ for any $x,y,z\in \fg$. Observe that the squaring is not involved here. 


We call the Lie superalgebra $\fg:=(\fg, [\cdot ,\cdot]_\fg,s_\fg)$ a~\textit{NIS-superalgebra} (sometimes called quadratics,  see \cite{BeB}) if it has a~non-degenerate invariant $\ev$-antisymmetric bilinear form ${\mathscr B}_\fg$.  A NIS-Lie superalgebra is called {\it irreducible} if it cannot be decomposed into a direct sum of mutually orthogonal ideals. 


For any Lie superalgebra $\fg$ in characteristic $2$, its \textit{derived algebras} are
defined to be (for $i\geq 0$)
\[
\fg^{(0)}: =\fg, \quad
\fg^{(i+1)}=[\fg^{(i)},\fg^{(i)}]_\fg+\Span\{s_\fg(x)\mid x\in (\fg^{(i)})_\od\}.
\]
A linear map ${\mathscr D}:\fg\rightarrow \fg$ is called a~\textit{derivation} of the Lie superalgebra $\fg$ if, in addition to
\begin{equation}
\label{Der1} {\mathscr D}([x,y]_\fg)=[{\mathscr D}(x),y]_\fg+[x,{\mathscr D}(y)]_\fg\quad \text{for any $x\in \fg_\ev$ and $y\in \fg$},
\end{equation}
we have 
\begin{equation}
\label{Der2} {\mathscr D}(s_\fg(x))=[{\mathscr D}(x),x]_\fg\quad \text{for any $x\in \fg_\od$}.
\end{equation}
It is worth noticing that condition (\ref{Der2}) implies condition (\ref{Der1}) if $x,y\in \fg_\od$. 

We denote the space of all derivations of $\fg$
by $\fder (\fg)$. 
\subsection{Non-degenerate $\od$-antisymmetric bilinear forms on Lie superalgebras for $p= 2$}

Let $\fg=\fg_\ev\oplus \fg_\od$ be a~Lie superalgebra in characteristic $2$.

Let $\omega:\fg \times \fg \rightarrow \Kee$ be a bilinear form.  We say that $\omega$ is {\it closed} if the following cocycle conditions are satisfied: 
\begin{eqnarray}
\label{Axiom2b} \omega ([x,y],z)+ \omega ([z,x],y)+ \omega ([y,z],x) & = & 0 \quad  \text{ for all $x,y,z\in \fg$}; \\
\label{Axiom3b} \omega(s_\fg(x),z) & = & \omega(x,[x,z])\quad  \text{ for all  $x\in \fg_\od$ and $z\in \fg$}.
\end{eqnarray}
Observe that the squaring has to be taken into account in defining the cocycle condition.  For more details,  see \cite{BGL2}. \\

\begin{mProposition}[Link NIS's with $\od$-antisymmetric forms]\label{NISvsSymp}
Let $(\fg, \mathscr B)$ be a~NIS Lie super-algebra in characteristic 2 such that $\mathscr B$ is homogeneous. Then $\fg$ admits a homogeneous,  non-degenerate, closed,  $\od$-antisymmetric  form $\omega$ if and only if there exists an invertible $\Delta\in \fder(\fg)$, where $p(\omega)=p(\mathscr B)+p(\Delta)$,  such that ${\mathscr B}$ is $\Delta$-invariant and 
\[
\omega(x,y)=\mathscr B (\Delta(x),y) \quad \text{for all $x,y\in \fg$}.
\]

\end{mProposition}
\begin{proof}

Let us first prove the ``only if'' part. Suppose that $\fg$ admits a non-degenerate,  homogeneous, closed, $\od$-antisymmetric form $\omega$. Since $\omega$ and ${\mathscr B}$ are non-degenerate,  there exists $\Delta\in \End(\fg)$ such that
\[
\omega(x,y)={\mathscr B} (\Delta(x),y) \quad \text{for all $x,y\in \fg$}.
\]
Let us prove that $\Delta$ is an invertible derivation. Indeed, the fact that $\Delta$ is invertible follows from the fact that $\omega$ is non-degenerate. Moreover, the fact that ${\mathscr B}$ is $\Delta$-invariant follows from condition (\ref{Axiom1}) or (\ref{AxiomO1b}) above.  For all $x,y,z\in \fg$,  we have
\begin{align*}
& {\mathscr B} (\Delta([x,y]),z) =  \omega([x,y],z)= \omega([z,x],y)+\omega([y,z],x) =  {\mathscr B} (\Delta([z,x]),y)+ {\mathscr B} (\Delta([y,z]),x)\\[2mm]
& =  {\mathscr B} ([z,x],\Delta(y))+ {\mathscr B} ([y,z],\Delta(x)) =  {\mathscr B} (z,[x,\Delta(y)])+ {\mathscr B} (z,[y,\Delta(x)])\\[2mm]
&= {\mathscr B} (z,[x,\Delta(y)]+ [y,\Delta(x)]).
\end{align*}
Therefore, $\Delta([x,y])=[\Delta(x),y]+[x, \Delta(y)] $ since ${\mathscr B}$ is non-degenerate and $\ev$-antisymmetric.  On the other hand,  (for all $x\in \fg_\od$ and $z\in \fg$):
\[
\begin{array}{l}
{\mathscr B} (\Delta(s(x)),z) =  \omega(s(x),z) =  \omega(x,[x,z])= {\mathscr B} (\Delta(x),[x,z]) = {\mathscr B} ([\Delta(x),x],z).
\end{array}
\]
It follows that $\Delta(s(x))=[x,\Delta(x)]$ since ${\mathscr B}$ is non-degenerate. Therefore, $\Delta$ is indeed a~derivation.  

For the ``if'' part,  the proof results from the fact that $\Delta$ is a derivation and the bilinear form ${\mathscr B}$ is $\Delta$-invariant.  \end{proof}


\subsection{Left-symmetric and left-alternative algebras in characteristic $2$}
In this subsection,  we will show that every homogeneous, non-degenerate, closed,  $\od$-antisymmetric form on a~Lie superalgebra $\fg$ in characteristic 2 naturally induces a new algebra structure, denoted by $\star$,  on the underlying vector space such that its commutator is nothing but the bracket on $\fg$.  Moreover, this structure is left-alternative, and therefore left-symmetric (see, Prop. \ref{leftalt}). These structures are also known as {\it Vinberg algebras} as they were used in \cite{V} (for a~thorough study of these structures,  see \cite{N}.) In addition, if $\fg$ has a NIS $\mathscr B$ then $\fg$ has an invertible derivation, relating  $\omega$ and $\mathscr B$,  such that the derivation will define the structure $\star$.  In Section  \ref{symplectic},  we will show that this derivation yields a Lie bisuperalgebra structure on $\fg$, and hence a Manin-triple $(\fg\oplus\fg^*, \fg, \fg^*)$.

Let $(\fg, \omega)$ be a~Lie superalgebra in characteristic 2 equipped with a homogeneous, non-degenerate,  closed, $\od$-antisymmetric form~$\omega$.

For all $x,y\in \fg$, the map $z\mapsto \omega (y,[x,z])$ is a~linear form on $\fg$. Moreover, the map $z\mapsto \omega(z,\cdot)$ identifies $\fg$ with $\fg^*$. It follows that there exists $x\star y \in \fg $ such that 
\[
\omega(x\star y,z)=\omega (y, [x,z]).
\]
Observe that if both $x$ and $y$ are homogeneous, then $x\star y$ is homogeneous as well.  Moreover, $p(x\star y)=p(x)+p(y)$. 

For all $x,y,z\in \fg$, we have
\[
0= \omega(x,[y,z])+\omega(y,[z,x])+\omega(z,[x,y])= \omega(y\star x,z)+\omega(x\star y,z)+\omega(z,[x,y]).
\]
Therefore, 
\begin{equation}
\label{symbra}
[x,y]=x\star y+ y \star x\;\text{ for all $x,y\in \fg$}.
\end{equation}
Aditionnally,  for all $x\in \fg_\od$ and $y\in \fg$, we have
\[
0= \omega(s(x),y)+\omega(x,[x,y])= \omega(x\star x,y).
\]
Therefore, 
\begin{equation}
\label{symsq}
s(x)=x \star x\;\text{ for all $x\in \fg_\od$}.
\end{equation}
Now, for all $x,y,z,t\in \fg$, we have
\[
\begin{array}{lcl}
0 &= & \omega(t, [x,[y,z]]+[y,[z,x]]+ [z,[x,y]])\\[2mm]
&= & \omega(t,[x,[y,z]]) +\omega(t,[y,[z,x]])+\omega(t, [z,[x,y]])\\[2mm]
&= & \omega(y \star ( x \star t),z) +\omega(x \star (y \star t)t,z)+\omega([x,y]\star t, z).
\end{array}
\]
Now Eq. (\ref{symbra}) and the fact that $\omega$ is non-degenerate imply that 
\[
\mathrm{Asso}(x,y,t)=\mathrm{Asso}(y,x,t)\; \text{ for all $x,y,t\in \fg$},
\]
where 
\[
\begin{array}{lcl}
\mathrm{Asso}(x,y,t)& := & (x\star y) \star t+ x \star (y \star t) \; \text{ for all $x,y,t\in \fg$}\\[2mm]
\end{array}
\]
On the other hand, 
\[
\begin{array}{lcl}
0 &= & \omega([s(x),y]+[x,[x,y]],t)\\[2mm]
&= & \omega(t,[s(x),y]) +\omega(t,[x,[x,y]])\\[2mm]
&= & \omega(s(x) \star t,y) +\omega(x \star ( x \star t ) ,y).
\end{array}
\]
Now Eq. (\ref{symsq}) and the fact that $\omega$ is non-degenerate imply that 
\[
\mathrm{Asso}(x,x,t)=0 \; \text{ for all $x\in \fg_\od$ and $t\in \fg$}.
\]
Suppose further that $\fg$ has a NIS $\mathscr B$ that  is homogeneous,  and let us describe $\star$ in terms of the derivation $\Delta$.  For all $x,y,z\in \fg$,  and using Proposition \ref{NISvsSymp}  we have
\[
0=\omega(x\star y ,z)-\omega (y,[x,z])={\mathscr B}(\Delta(x\star y),z)-{\mathscr B}(\Delta(y),[x,z])={\mathscr B}(\Delta(x\star y)-[\Delta(y),x],z).
\]
Since ${\mathscr B}$ is non-degenerate and $\Delta$ is invertible, it follows that
\[
x\star y=\Delta^{-1}([x,\Delta(y)])\; \text{ for all $x,y\in \fg$}.
\]
In particular, $s(x)=x\star x= \Delta^{-1}([x,\Delta(x)])$. We then arrive at the following definition.\\

An {\bf algebra} $\fg$ in characteristic $2$ is called {\it left-symmetric} (resp. {\it left-alternative}) if it is endowed with a~bilinear operation $\star$ satisfying (for all $x,y,z\in \fg$):
\[
\mathrm{Asso}(x,y,z) =  \mathrm{Asso}(y,x,z)\quad \text{(resp. $\mathrm{Asso}(x,x,z)=0$)},
\]
where
\[
\begin{array}{lcl}
\mathrm{Asso}(x,y,z)& := & (x\star y) \star z+ x \star (y \star z) \text{ for all $x,y,z\in \fg$}, 
\end{array}
\]
is called the {\it associator}.\\

\begin{mProposition}[left-alternative implies left-symmetric] \label{leftalt}
In characteristic $2$, every left-alternative algebra is left-symmetric.

\end{mProposition}
\begin{proof}

Indeed, for all $x,y,z\in \fg$, we have
\[
\begin{array}{l}
\mathrm{Asso}(x+ y, x + y,z)=  \mathrm{Asso}(x, x,z)+ \mathrm{Asso}(x, y,z)+ \mathrm{Asso}(y, x,z)+ \mathrm{Asso}(y, y,z). \qed
\end{array}
\]
\noqed\end{proof}
\begin{mRemark}{\rm (i) This proposition shows that left-alternatives algebras constitute a subclass of left-symmetric algebras, which is not true if $p\not =2$.  Whether these classes are the same or not is an open question. 

(ii) Note that the structure $\star$ associated with a closed $\od$-antisymmetric ortho-orthogonal (or periplectic) form is left-alternative.

}

\end{mRemark}

\begin{mTheorem}[Constructing Lie superalgebras from Lie algebras]\label{WNIS} Let $p=2$.  
\textup{(}i\textup{)} Every Lie algebra $\fg$ with an invertible derivation $\Delta\in \fder(\fg)$ admits a~left-symmetric structure. Explicitly, 
\begin{equation}
\label{deltastar}
x\star y :=\Delta^{-1}([x, \Delta(y)])\; \text{ for all $x,y\in \fg$}.
\end{equation}

\textup{(}ii\textup{)}  Every $\Zee/2$-graded Lie algebra $\fg:=\fg_\ev\oplus \fg_\od$ with an invertible homogeneous derivation $\Delta$  can be turned into a~Lie superalgebra provided that
\[
\mathrm{Asso}(x,x,y)=0\; \text{ for all $x\in \fg_\od$ and $y\in \fg$},
\] 
where the assciator $\mathrm{Asso}$ is associated with the operation \eqref{deltastar}. Explicitly, the map  
\[
s: \fg_\od \rightarrow \fg_\ev \quad x\mapsto x \star x
\]
is a~squaring on $\fg$.

\end{mTheorem}
\begin{proof}

(i) Let us first compute the associator. For all $x,y,z\in \fg$, we have
\begin{align*}
& \mathrm{Asso}(x,y,z)  =  (x\star y) \star z+ x \star ( y \star z) = \Delta^{-1}\left ([x\star y, \Delta(z)]+ [x, \Delta(y \star z)]\right )\\[2mm]
& =  \Delta^{-1}\left ([\Delta^{-1}[x, \Delta(y)], \Delta(z)]+ [x, [y, \Delta(z)]\right ) \\[2mm]
& =  \Delta^{-1}\left ([[x,y], \Delta(z)]+[\Delta^{-1}[\Delta(x),y],\Delta(z)]+[x, [y, \Delta(z)]]\right ).
\end{align*}
Now, using the Jacobi identity and the fact that $\Delta$ is a~derivation we see that
\[
\begin{array}{lcl}
\mathrm{Asso}(x,y,z) + \mathrm{Asso}(y,x,z)& = & \Delta^{-1}\left ([\Delta^{-1}[\Delta(x),y],\Delta(z)]+[x, [y, \Delta(z)]]\right )\\[2mm]
&&+\Delta^{-1}\left ([\Delta^{-1}[\Delta(y),x],\Delta(z)]+[y, [x, \Delta(z)]]\right )\\[2mm]
&=&\Delta^{-1}\left ([[x,y],\Delta(z)]+[\Delta(z),[x,y]]        \right ) =0.
\end{array}
\]
Therefore, the operation (\ref{deltastar}) defines a~left-symmetric structure on $\fg$.

(ii) Suppose now that $\fg=\fg_\ev\oplus \fg_\od$ is $\Zee/2$-graded and that $\Delta$ is homogeneous.  Let us check that the map
\[
s: \fg_\od \rightarrow \fg_\ev \quad x\mapsto x \star x
\]
is indeed a~squaring. For all $\lambda \in \Kee$ and $x\in \fg_\od$, we have
\[
s(\lambda x)=(\lambda x)\star (\lambda x)= \lambda^2 (x\star x)=\lambda^2 s(x).
\]
Additionally,  for all $x,y \in \fg_\od$:
\begin{align*}
& s(x+y)-s(x)-s(y)=(x+y)\star (x+y)-x\star x -y \star y = x\star y + y\star x\\[2mm]
&=\Delta^{-1}[x, \Delta(y)]+ \Delta^{-1}[y,\Delta(x)]=\Delta^{-1}(\Delta([x,y]))= [x,y].
\end{align*}
Moreover, for all $x\in \fg_\od$ and $y\in \fg$, we have
\[
\begin{array}{lcl}
[s(x), y] + [x,[x,y]] & = & s(x)\star y + y \star s(x) +  x\star [x,y] + [x,y]\star x\\[2mm]
& = & (x \star x)\star y+ y \star (x\star x)+ x \star (x\star y+ y \star x)+ (x\star y+ y \star x) \star x\\[2mm]
& = & (x \star x)\star y+ x \star (x\star y)+ y \star (x\star x) + (y \star x )\star x+ x \star (y \star x) \\[2mm]
&&+ (x\star y) \star x\\[2mm]
&=& \mathrm{Asso}(x,x,y) + \mathrm{Asso}(y,x,x) +  \mathrm{Asso}(x,y,x)= \mathrm{Asso}(x,x,y).
\end{array}
\]
It follows that the Jacobi identity involving the squaring holds if and only if 
\[ 
\mathrm{Asso}(x,x,y)=0 \text{ for all $x\in \fg_\od$ and $y\in \fg$}. \qed
\]
\noqed
\end{proof}

\begin{mCorollary}[The case of $2$-step nilpotent algebras] Every $\Zee/2$-graded $2$-step nilpotent Lie algebra in characteristic $2$ with an invertible  homogeneous derivation can naturally be turned into a~Lie superalgebra.

\end{mCorollary}
\begin{proof}

Let us first observe that since $\Delta$ is invertible, then $\Delta([\fg,\fg])=[\fg,\fg]$. Now, for all $x\in \fg_\od$ and  $y\in \fg$:
\begin{align*}
& \mathrm{Asso}(x,x,y)  =  x\star (x \star y)+(x \star x)\star y = x\star \Delta^{-1}[x, \Delta(y)] + \Delta^{-1}[x,\Delta(x)]\star y\\[2mm]
&=  \Delta^{-1}\left ([x, [x, \Delta(y)]] + [\Delta^{-1}[x,\Delta(x)], \Delta(y)]\right )= \Delta^{-1}(0)=0.\qed
\end{align*}
\noqed
\end{proof}

\section{The Manin triples and matched pairs for Lie superalgebras in characteristic $2$}
Throughout this section $\mathbb K$ is an arbitray field of characteristic $2$.
\subsection{Manin triples }
Let $\fh$ be a~Lie superalgebra and $\fg, \fk$ be two Lie sub-superalgebras of $\fh$. The triple $(\fh, \fg, \fk)$ is called a~\textit{Manin triple} if it satisfies the following two conditions:
\begin{enumerate}
\item $\fh =\fg\oplus \fk$,  the direct sum of vector subspaces $\fg$ and $\fk$;

\item there exists an~even NIS $\mathscr B$ on $\fh$ such that $\fg$ and $\fk$ are isotropic with respect to $\mathscr B$.
\end{enumerate}
For a~comprehensive construction of Manin triples for complex Lie algebras, see \cite{CP}. For Lie superalgebras in characteristic 2, the first example of a~Manin triple was given in \cite{BeBou}. Here we generalize the construction developed in \cite{BeBou}.\\

\begin{mTheorem}[On Manin Triples] \label{mainmanin}
Let $(\fg, [\cdot,\cdot]_\fg,s_\fg)$ be a~Lie superalgebra in characteristic $2$, and suppose that on the dual space $\fg^*$ a~Lie superalgebra structure $(\fg^*, [\cdot,\cdot]_{\fg^*},s_{\fg^*})$ is also given. On the space $\fh:=\fg\oplus \fg^*$, we consider a~natural even bilinear form defined by
\[
{\mathscr B}_\fh(x+f,y+g)=g(x)+f(y) \quad \text{ for every $x,y\in \fg$ and $f,g\in \fg^*$.}
\] 
There exists a~unique Lie NIS-superalgebra structure on $(\fh, {\mathscr B}_\fh)$  provided the conditions \eqref{Sq}, \eqref{Sq*}, and \eqref{Bra} are fulfilled:
\begin{eqnarray}\label{maninsq1}
\label{Sq} s_\fg(x)\circ \ad_h&=&[x,x\circ \ad_h]_\fg+x\circ \ad_{h \circ \ad_x}\; \text{ for all $x\in \fg_\od$ and $h\in \fg^*$};\\ 
\label{Sq*}
s_{\fg^*}(f)\circ \ad_y&=&[f,f\circ \ad_y]_{\fg^*}+f\circ \ad_{y \circ \ad_f}\; \text{ for all $y\in \fg$ and $f\in \fg^*_\od$};\\
\label{Bra}
[x\circ \ad_f, y]_\fg&=&y \circ \ad_{f\circ \ad_x}+[x,y\circ \ad_f]_\fg+x \circ \ad_{f\circ \ad_y}+[x,y]_\fg\circ \ad_f,
\end{eqnarray}
\begin{equation*}
\label{Cond1}
\begin{array}{l}
\text{for all $x,y \in \fg_\ev$, and $f\in \fg^*$}; \\ 
\text{or for all $x\in \fg_\ev$, $y\in \fg$ and $f\in \fg^*_\ev$};\\
\text{or for all $x\in \fg_\od$, $y\in \fg$ and $f\in \fg^*_\od$}.
\end{array}
\end{equation*} 
Explicitly,  
\begin{itemize}
\item The bracket is given by \textup{(}for all $x,y\in \fg$ and for all $f,h\in \fg^*$\textup{)}
\[
[x+f,y+h]_\fh:=[x,y]_\fg+y\circ \ad_f+x\circ \ad_h+[f,h]_{\fg^*}+f\circ \ad_y+ h\circ \ad_x.
\]
\item The squaring is given by
\[
s_\fh(x+f)=s_\fg(x)+s_{\fg^*}(f)+f\circ \ad_x+ x\circ \ad_f\; \text{ for all $x\in \fg_\od$ and $f\in \fg^*_\od$.}
\]
\end{itemize}
The bilinear form ${\mathscr B}_\fh$ is $\ad_\fh$-invariant and even.

\end{mTheorem}
Clearly,
\[
[\cdot,\cdot]_\fh|_\fg=[\cdot,\cdot]_\fg, \quad [\cdot,\cdot]_\fh|_{\fg^*}=[\cdot,\cdot]_{\fg^*}, \quad s_\fh |_\fg=s_\fg, \quad s_\fh |_{\fg^*}=s_{\fg^*},
\]
and therefore $(\fh:=\fg\oplus \fg^*,\fg,\fg^*)$ is a~Manin triple.

For the proof,  see Section \ref{mpairs}.

The reader may wonder why there is no analogue of Eq. \eqref{Bra} on $\fg^*$. The answer is yes indeed it should be,  but Lemma \ref{EqBra0} shows that it is equivalent to Eq. \eqref{Bra}.\\

\begin{mLemma}[Technical]\label{EqBra0} \textup{(}i\textup{)} For all $x,y\in \fg$ and $f,h\in \fg^*$, we have
\[
h([x,y\circ \ad_f]_\fg)=[f,h\circ \ad_x]_{\fg^*}(y).
\]
\textup{(}ii\textup{)} For all $x,y\in \fg$ and $f,h\in \fg^*$, the following two conditions are equivalent
\[
\begin{array}{lcl}
[x\circ \ad_f, y]_\fg+y \circ \ad_{f\circ \ad_x}&=&[x,y\circ \ad_f]_\fg+x \circ \ad_{f\circ \ad_y}+[x,y]_\fg\circ \ad_f,\\[2mm]
[f\circ \ad_x, h]_{\fg^*}+h \circ \ad_{x\circ \ad_f}&=&[f,h\circ \ad_x]_{\fg^*}+f \circ \ad_{x\circ \ad_h}+[f,h]_{\fg^*}\circ \ad_x.
\end{array}
\]

\end{mLemma}
\begin{proof} 

For Part (i), we have
\[
\begin{array}{rcl}
h([x,y\circ \ad_f]_\fg)&=& h\circ \ad_x (y\circ \ad_f)=y\circ \ad_f(h\circ \ad_x) \quad \text{(we identify $\fg$ with $\fg^{**}$)}\\[2mm]
&=&y([f,h\circ \ad_x]_{\fg^*})=[f,h\circ \ad_x]_{\fg^*}(y).
\end{array}
\]

Part (ii) follows by evaluating the first equation at an arbitrary $h\in \fg^*$, evaluating the second equation at an arbitrary $y\in \fg$,  and then using Part (i).\end{proof}

\subsection{A few examples of Manin triples in characteristic $2$}\label{examples}
Here we give three examples of Manin triples of Lie superalgebras in characteristic $2$. \\

\noindent {\bf Example 1:} A Manin triple constructed from the Heisenberg superalgebra.  Consider the Heisenberg superalgebra $\fg:=\fhei(0|2)$ spanned by  $p,q$ (odd) and $z$ (even), with $s_\fg(p)=s_\fg(q)=0$ and the only nonzero bracket: $[p,q]=z$. We endow $\fg^*$ with a~Lie superalgebra structure given by:
\[
[p^*, z^*]= s q^*+ t p^* , \quad [q^*, z^*]=u q^*+ v p^* \text{  for any fixed $s,t,u,v\in \mathbb{K}$}.
\]
A NIS-Lie superalgebra $\fh:=\fhei(0|2)\oplus \fhei(0|2)^*$ exists.  A direct computation shows the following:
\begin{itemize}
\item The bracket on $\fh$ is given by 
\[
\begin{array}{l}
[p,p^*]= t z, \quad [p,q^*]=v z, \quad [p,z^*]=t p+v q+q^*,\\[2mm]
[q,p^*]= s z, \quad [q,q^*]=u z, \quad [q,z^*]=s p+u q+p^*.
\end{array}
\]
\item The squaring on $\fh$ is given by (where $\alpha_1, \alpha_2, \alpha_3, \alpha_4\in \mathbb{K}$):
\[
\begin{array}{lcl}
s_\fh(\alpha_1 p+\alpha_2q+\alpha_3 p^*+\alpha_4 q^*)&=&(\alpha_1\alpha_2+\alpha_1\alpha_3 t+\alpha_1\alpha_4 v+\alpha_2\alpha_3s+\alpha_2\alpha_4u)z.\\[2mm]
\end{array}
\]
\end{itemize}
We give an example of an outer derivation ${\mathscr D}$ on $\fh$ that is invertible. This example will justify the studies in Section \ref{symplectic}. The example we provide here was borrowed from \cite{BeBou} in the case where $ \fhei(0|2)^*$ has an abelian structure; namely, $s=t=u=v=0$. The ground field $\Kee$ is assumed to be algebraically closed, hence infinite. The invertible derivation is as follows (for $\alpha_9, \alpha_{10}\not=0,1$):
\[
{\mathscr D}=\alpha_{2} D_2+ \alpha_{4} D_4+ \alpha_9 D_9+ \alpha_{10}D_{10}+ (\alpha_9+ \alpha_{10})D_{11}, 
\]
where 
\[
\begin{array}{lcllcllcllcl}
D_2&=&q^*\otimes  \widehat{q}, & D_4&=&p^*\otimes  \widehat{p},\\[2mm]
D_9&=&p\otimes  \widehat{p}+q^*\otimes  \widehat{q^*}+z\otimes  \widehat{z}, & D_{10}&=&q\otimes  \widehat{q}+p^*\otimes  \widehat{p^*}+z\otimes  \widehat{z},\\[2mm]
 D_{11}&=&q^*\otimes  \widehat{q^*}+p^*\otimes  \widehat{p^*}+z^*\otimes  \widehat{z^*}.
\end{array}
\]
{\bf Example 2:} The Poisson superalgebras \cite{LSh}.  Consider $\fg = \fpo(0|2n - 1)$. Let $\fh = \fpo(0|2n)$ be generated by $\Bbb K[\theta_1,...,\theta_{2n}]$ and $\fg$ by $\Bbb K[\theta_1,...,\theta_{2n-1}]$. The form ${\mathscr B}_\fh$ on $\fh$ is given by the Berezin integral ${\mathscr B}_\fh(x,y) =\int xy\; \mathrm{vol}(\theta)$. Let $\fg^* = \fg \cdot \zeta$, where $\zeta = \theta_{2 n} +   \sum_{i=1}^{2n-1} k_i \theta_i$ and $\sum_{i=1}^{2n-1}{k _i}^2 = 1.$

\noindent {\bf Example 3:} An example of a Lie superalgebra endowed with both NIS and a closed $\od$-antisymmetric ortho-orthogonal form.  Let $(\fh, [\cdot,\cdot]_\fh, s_\fh)$ be a Lie superalgebra with an invertible derivation~$\mathscr D$ (take for instance $\fh=\fhei(0|2)\oplus \fhei(0|2)^*$ from Example 1).   Consider the trivial structure on the dual space $\fh^*$.  There exists a NIS $\mathscr B$ on the Lie superalgebra $\fg:=\fh\oplus \fh^*$ given as in Theorem \ref{mainmanin}.  Let us extend the derivation $\mathscr D$ to an invertible derivation $\Delta$ on  $\fg$.  Indeed, let us define
\[
\Delta (x+f):={\mathscr D}(x) + f \circ {\mathscr D} \quad  \text{ for all $x\in \fh$ and $f\in \fh^*$}.
\] 
 Moreover,  $\Delta$ is $\mathscr B$-symmetric. Therefore, the form $\omega(\cdot,\cdot)={\mathscr B}(\Delta(\cdot), \cdot)$ is closed, $\od$-antisymmetric and ortho-orthogonal on $\fg$ by Prop. \ref{NISvsSymp}.   It is worth mentioning that this example supports the study in Section~\ref{symplectic}. 
\subsection{Review on cohomology of Lie superalgebras in characteristic $2$} 
Recall basic facts concerning the cohomology of Lie superalgebras in characteristic $2$; for more details, see \cite{BGL2}. 1-cochains are just linear functions on $\fg$ with values in a~$\fg$-module $M$. In our context, the {\it differential} of a~1-cochain $c$ is defined as follows:
\begin{equation}
\label{diffsuper}
\begin{array}{lcl}
 {\mathfrak d}c(x,y) & = & c([x,y]) + x \cdot c(y) + y \cdot c(x)\quad \text{for all $x,y\in \fg$};\\[1mm]
{\mathfrak d}c(x,x) & = & c(s_\fg(x))+x\cdot c(x) \quad \text{for all $x\in \fg_\od$}.
\end{array}
\end{equation}
A 1-cocycle $c$ on $\fg$ with values in a~$\fg$-module $M$ must satisfy the following conditions:
\begin{eqnarray}
\label{Cond1} x\cdot c(y)+y\cdot c(x)+c([x,y]_\fg)&=& 0 \quad \text{for all $x,y\in \fg$},\\[2mm]
\label{Cond2} x\cdot c(x)+c(s_\fg(x)) &=& 0\quad \text{for all $x\in \fg_\od$}.
\end{eqnarray}
The space of all 1-cocycles is denoted by $Z^1(\fg;M)$.

To every Lie superalgebra $\fg$ in characteristic $2$, we can assign its {\it desuperization} denoted by $F(\fg)$, where $F$ is the functor of forgetting the super structure. The space $F(\fg)$ turns $\fg$ into an ordinary $\mathbb Z/2$-graded Lie algebra; for more details, see \cite{BLLSq}.  Now, any 1-cocycle  $c$ on the Lie algebra $F(\fg)$ is captured only by the condition \eqref{Cond1},  and therefore the cohomology of $F(\fg)$ viewed as a~Lie  algebra and that of $\fg$ viewed as a~Lie superalgebra are not necessarily isomorphic. We refer to \cite{BGL2} for illustrative examples.

Let us define a~bilinear form $\langle - ,  - \rangle: \fg^*\times \fg \rightarrow \mathbb K$ by $\langle f , x\rangle=f(x)$ for all $x\in \fg$ and for all $f\in \fg^*$. Since $\fg$ is finite-dimensional,  we have $(\fg \otimes \fg)^*= \fg^*\otimes \fg^*$ so we also denote by $\langle -,  - \rangle$ the bilinear form on $(\fg \otimes \fg)^*\times (\fg\otimes \fg)$ defined by 
\[
\langle f\otimes h, x\otimes y\rangle=f(x) h(y) \text{ for all $x,y\in \fg$ and $f,h\in \fg$},
\]
then extend  the definition  by linearity. Given the bracket $[\cdot, \cdot]_{\fg^*}:\fg^*\otimes \fg^* \rightarrow \fg^*$, its dual is the map $c_\fg:\fg\rightarrow \fg \otimes \fg$. Similarly, given the bracket $[\cdot, \cdot]_\fg:\fg\otimes \fg \rightarrow \fg$, its dual is the map $c_{\fg^*}:\fg^*\rightarrow \fg^* \otimes \fg^*$. More precisely, 
\begin{eqnarray}
\label{brag} \langle [f,h]_{\fg^*}, x\rangle & = & \langle f\otimes h, c_\fg(x)\rangle\quad  \text{ for all $x\in \fg$ and $f,h\in \fg^*$},\\[2mm]
\label{bragstar}\langle c_{\fg^*}(f), x\otimes y \rangle & =& \langle f, [x,y]_\fg\rangle \quad \text{ for all $x,y\in \fg$ and $f\in \fg^*$}.
\end{eqnarray}
We observe here that both $c_\fg$ and $c_{\fg^*}$ have images generated by symmetric  2-tensors. \\


\begin{mProposition}[A cohomological interpretation]\label{cohinterp}
The conditions \eqref{Sq}, \eqref{Sq*} and \eqref{Bra} in Theorem $\ref{mainmanin}$ are cohomological. Namely, they are equivalent to the fact that $c_\fg\in Z^1(\fg; \fg\otimes \fg)$ and $c_{\fg^*} \in Z^1(\fg^*; \fg^*\otimes \fg^*)$.

\end{mProposition} 
\begin{proof} 

We will only be dealing with the map $c_\fg$.  Studying the properties of the map $c_{\fg^*}$ is absolutely the same. 

Now, checking that the cocycle condition \eqref{Cond1} is equivalent to \eqref{Bra} is a~routine, so we omit this part. Let us check the conditions involving the squaring. We have
\begin{align*}
& \langle f\otimes h, x \cdot c_\fg(x)+c_\fg(s_\fg(x))\rangle=\langle f\otimes h, x \cdot c_\fg(x)\rangle+\langle f\otimes h, c_\fg(s_\fg(x))\rangle \\[2mm]
& = x\cdot (f \otimes h)(c_\fg(x))+(f\otimes h)(c_\fg(s_\fg(x)))= ((x\cdot f)\otimes h +f \otimes (x\cdot h))(c_\fg(x))+[f,h]_{\fg^*}(s_\fg(x)) \\[2mm]
& =\langle [x\cdot f, h]_{\fg^*}, x\rangle+\langle [f,x\cdot h]_{\fg^*}, x \rangle+[f,h]_{\fg^*}(s_\fg(x))\\[2mm]
&= \langle [f\circ \ad_x,h]_{\fg^*}, x \rangle +\langle f,h\circ \ad_x]_{\fg^*}, x \rangle +\langle [f,h]_{\fg^*}, s_\fg(x)\rangle.
\end{align*}
On the other hand, let us consider Eq. \eqref{Sq} applied to $f$. We get:
\begin{align*}
& \left ( s_\fg(x)\circ \ad_h+[x,x\circ \ad_h]_\fg+x\circ \ad_{h\circ \ad_x} \right ) (f)\\
&= [h,f]_{\fg^*}(s_\fg(x))+ f([x,x\circ \ad_h]_\fg)+\ad_{h\circ \ad_x}(f)(x)\quad  \text{ using Lemma (\ref{EqBra0})}\\
&=\langle [h,f]_{\fg^*}, s_\fg(x) \rangle+\langle [h, f \circ \ad_x]_{\fg^*}, x \rangle+\langle [h\circ \ad_x, f]_{\fg^*}, x\rangle.
\end{align*}
It follows that Eq. \eqref{Sq} is equivalent to Eq. \eqref{Cond2}.\end{proof} 
\begin{mProposition}[Relation between $c_\fg$ and $c_{\fg^*}$] \label{equicg} For all $f,g\in \fg^*$ and $x,y\in \fg$, we have 
\[
\langle f\cdot c_{\fg^*}(h)+h\cdot c_{\fg^*}(f)+c_{\fg^*}([f,h]), x \otimes y\rangle = \langle f\otimes h, x\cdot c_\fg(y)+y \cdot c_\fg(x)+c_\fg([x,y])\rangle.
\]
Thus, the map $c_\fg$ satisfies condition \eqref{Cond1} if and only if the map $c_{\fg^*}$ satisfies condition~\eqref{Cond1}.

\end{mProposition}

\begin{Remark} For $p\not =2$, the 1-cocycle conditions of the maps $c_\fg$ and $c_{\fg^*}$ are actually equivalent thanks to an analogue of Prop. \ref{equicg}.  For $p=2$,  Prop. \ref{equicg} is not enough to conclude that they are equivalent, due to the presence of the squaring in Eq. (\ref{Cond2}); conjecturally, they are not.
\end{Remark}

In order to prepare for the proof of Proposition \ref{equicg}, we will be using {\it Sweedler's notations}. Let us write 
\[
c_{\fg^*}(f):=\sum_{(f)} f_{(1)} \otimes f_{(2)} \text{ for all $f\in \fg^*$ and $c_{\fg}(x) :=\sum_{(x)} x_{(1)} \otimes x_{(2)}$ for all $x\in \fg$}.
\]
\begin{mLemma}[Technical] \label{equicgl} For all $x,y\in \fg$ and $f,h\in \fg^*$, we have
\[
\langle f \cdot c_{\fg^*}(h), x\otimes y\rangle= \langle f \otimes h, \sum_{(x)} x_{(1)} \otimes [x_{(2)},y] +\sum_{(y)}  y_{(1)} \otimes [x, y_{(2)}] \rangle.
\]

\end{mLemma}
\begin{proof} 

For all $f,h\in \fg^*$ and $x,y\in \fg$, we have 
\begin{align*}
&\langle f \cdot c_{\fg^*}(h), x\otimes y \rangle =\langle \sum_{(h)} ( [f, h_{(1)}]_{\fg^*} \otimes h_{(2)}+\ h_{(1)} \otimes [f,  h_{(2)}]_{\fg^*}, x\otimes y \rangle \\[1mm]
&=\sum_{(h)} ( f \otimes h_{(1)}) (c_\fg(x))\, h_{(2)}(y)+\ h_{(1)} (x) \, (f \otimes  h_{(2)})(c_\fg(y))\\[1mm]
&=\sum_{(h)} \sum_{(x)}  f(x_{(1)}) h_{(1)}(x_{(2)}) \, h_{(2)}(y)+\sum_{(h)} \sum_{(y)} h_{(1)} (x) \, f(y_{(1)}) h_{(2)}(y_{(2)})\\[1mm]
&=f( \sum_{(x)} c_{\fg^*}(h) (x_{(2)} \otimes y) \,x_{(1)} )+f( \sum_{(y)} c_{\fg^*}(h)( x \otimes y_{(2)})\, y_{(1)})\\[1mm]
&=f( \sum_{(x)} h([x_{(2)}, y]) \,x_{(1)} )+ \sum_{(y)} h([x, y_{(2)}])\, y_{(1)}])\\[1mm]
&=\langle f\otimes h, \sum_{(x)} x_{(1)}  \otimes [x_{(2)}, y]+\sum_{(y)}  y_{(1)} \otimes [x, y_{(2)}] \rangle. \qed
\end{align*}
\noqed
\end{proof}

\begin{proof}[(Proof of Proposition $\ref{equicg}$)] Using Lemma \ref{equicgl} and the fact that $c_\fg$ takes values in the space of symmetric 2-tensors, we have
\begin{align*}
&\langle f \cdot c_{\fg^*}(h)+ h\cdot c_{\fg^*}(f)+c_{\fg^*}([f,h]), x\otimes y \rangle \\[1mm]
&=(f\otimes h)( \sum_{(x)} x_{(1)}  \otimes [x_{(2)}, y]+\sum_{(y)}  y_{(1)} \otimes [x, y_{(2)}])\\[1mm]
&+ (f \otimes h) ( \sum_{(x)}  [x_{(2)}, y] \otimes x_{(1)}+  \sum_{(y)}  [ x, y_{(2)}]\otimes y_{(1)} )+(f\otimes h)(c_\fg([x,y])) \\[1mm]
&=(f\otimes h)( \sum_{(x)} x_{(1)}  \otimes [x_{(2)}, y]+ \sum_{(x)}  [x_{(2)}, y] \otimes x_{(1)}   \\[1mm]
&+ (f \otimes h) (\sum_{(y)}  y_{(1)} \otimes [x, y_{(2)}] +  \sum_{(y)}  [ x, y_{(2)}]\otimes y_{(1)} )+(f\otimes h)(c_\fg([x,y])) \\[1mm]
&=(f\otimes h)( \sum_{(x)}   [y, x_{(1)}]\otimes x_{(2)}+ \sum_{(x)}  x_{(1)}\otimes [x_{(2)}, y]  \\[1mm]
&+ (f \otimes h) (\sum_{(y)}  [x, y_{(1)}] \otimes y_{(2)} +  \sum_{(y)}  y_{(1)}  \otimes [ x, y_{(2)}])+(f\otimes h)(c_\fg([x,y]))\\[1mm]
& =(f\otimes h) (x\cdot c_\fg(y)+y \cdot c_\fg(x)+c_\fg([x,y])). \qed
\end{align*}
\noqed
\end{proof}

\subsection{$r$-matrices and Manin triples} \label{rmatManin}
This section is devoted to the study of triangular Lie bi-superalgebras in characteristic $2$. We refer to \cite{CP, D, GZB, M2} in the case of Lie (super)algebras defined over a~field of characteristic zero, and to \cite{F1, F2} in the case of modular Lie algebras.

A {\it Lie coalgebra} over a~field $\mathbb{K}$ of characteristic $p\not =2$ is a~vector space $\fg$ over $\mathbb{K}$ together with a~linear mapping $c: \fg \rightarrow \fg \otimes \fg$, such that
\begin{itemize}
\item[(i)] $\IM(c)\subseteq \IM( 1 - \tau)$, and
\item[(ii)] $(1 + \xi + \xi^2) \circ  (1 \otimes c) \otimes c =0,$ 
where $\tau : \fg \otimes \fg \rightarrow  \fg\otimes \fg$ denotes the switch mapping sending $x \otimes y$ to $y\otimes x$ for all
$x,y\in \fg$  while $\xi: \fg\otimes \fg \otimes \fg \rightarrow \fg\otimes \fg \otimes \fg$ denotes the cycle mapping sending $x\otimes y \otimes z$ to
$y \otimes z \otimes x$ for all $x, y,z \in \fg$. 
\end{itemize}The mapping $c$ is called the cobracket of $\fg$, (i) is called
{\it co-anticommutativity}, and (ii) is called the {\it co-Jacobi identity}.

For $p=2$, the co-anticommutativity condition (i) has to be changed to the following condition: 
\[
c^t(f\otimes f)=0 \text{ for all $f\in \fg^*$},
\]
where $c^t: \fg^*\otimes \fg^* \rightarrow \fg^*$ is the dual of $c$. Our choice is based on the fact that a~coalgebra structure on $\fg$ will define a~Lie algebra structure on the dual space $\fg^*$.  Recall that for $p=2$, the anti-commutativity of the bracket is expressed as $[a,a]=0.$ 

A $\mathbb{Z}\slash 2 $-graded Lie coalgebra is a~$\mathbb{Z}\slash 2 $-graded vector space  $\fg=\fg_\ev\oplus \fg_\od$ over $\mathbb K$ endowed with a~cobracket $c$ such that $c(\fg_{\bar i})\subset (\fg\otimes \fg)_{\bar i}$, where $\bar i=\ev, \od$ (namely, $c$ is even). 

A {\it Lie bialgebra} over $\mathbb{K}$ is a~vector space $\fg$ over $\mathbb{K}$ together with $\mathbb{K}$-linear mappings
$[\cdot,\cdot]: \fg \otimes \fg \rightarrow \fg$ and $c: \fg \rightarrow \fg\otimes \fg$ such that $(\fg,[\cdot, \cdot])$ is a~Lie algebra, $(\fg, c)$ is a~Lie
coalgebra and $c$ is a~derivation of the Lie algebra $\fg$ to the $\fg$-module $\fg \otimes \fg$. Namely, 
\[
c([x,y])=x\cdot c(y) - y\cdot c(x) \text{ for all $x,y\in \fg$}.
\]

A {\it coboundary} Lie bialgebra  is a~Lie bialgebra such that the cobracket $c$ 
is a cboundary, i.e., there exists an element $r\in \fg\otimes \fg$ such that
$c(x)=x \cdot r$ for all $x\in \fg$.  

Let $\fg$ be a~Lie algebra and $r \in  \fg \otimes \fg$ be defined by $r =\sum_{i=1}^n a_i \otimes  b_i$, with $a_i \in  \fg,$ and $b_i \in  \fg$ for all $i \in \{1, \ldots ,n\}$. Consider $\fg$ embedded in its universal 
enveloping superalgebra $U(\fg)$. We put 
\[
r^{12} =\sum_{1\leq i \leq n} a_i \otimes  b_i \otimes 1, \quad r^{13} =\sum_{1\leq i \leq n} a_i \otimes 1 \otimes  b_i, \quad \text{ and } \quad r^{23} =\sum_{1\leq i \leq n} 1 \otimes  a_i \otimes  b_i.
\] 
The mapping $\mathrm{CBY}:\fg \otimes \fg \rightarrow \fg\otimes \fg \otimes \fg$ given by 
\[
r \rightarrow [r^{12},r^{13}]+[r^{12},r^{23}]+[r^{13},r^{23}]
\] 
is called the {\it classical Yang-Baxter map }, the equation $\mathrm{CYB}(r)=0$ is called the {\it classical Yang-Baxter equation} (CYBE), and a~solution of the CYBE is called a~{\it classical}
$r$-matrix.

If $\fg$ is a~Lie algebra and $V$ is a~$\fg$-module, then the set of $\fg$-invariant elements of $V$ is defined as
\[
V^\fg:=\{v \in V\, |\, x\cdot v=0 \text{ for all }  x\in \fg\}.
\]
In the case of Lie algebras defined over a~field of characteristic zero,  Drinfeld proved the following theorem.\\

\begin{mTheorem}[Lie bialgebra structure,  see\cite{D}] \label{D} Let $r\in \fg\otimes \fg$ and define $\delta_r (x) := x \cdot r$ for every $x \in  \fg$. Then, $r$ defines a~Lie bialgebra structure on $\fg$ if and only if 

\textup{(}i\textup{)} $r+\tau(r) \in (\fg \otimes \fg)^\fg$, and

\textup{(}ii\textup{)} $\mathrm{CYB}(r) \in (\fg \otimes \fg \otimes \fg)^\fg$.
\end{mTheorem} 

Therefore, every solution $r$ of the CYBE satisfying $r+\tau(r) \in (\fg \otimes \fg)^\fg$ gives rise to a~coboundary
Lie bialgebra structure on $\fg$ by means of $\delta_r$. Following \cite{D}, such a~Lie bialgebra structure is called {\it quasi-triangular}. The quasi-triangular Lie bialgebra structures arising from anti-symmetric
classical $r$-matrices are called {\it triangular}.

Following \cite{F1, F2}, coboundary Lie bialgebra structures defined over a~field of characteristic $2$ are called {\it quasi-triangular} if $\mathrm{Im}(\delta_r) \subseteq \mathrm{Im}( 1 + \tau)$, and called {\it triangular} if $r\in \mathrm{Im}(1-\tau)$.

For $p=2$, the condition (i) in Theorem \ref{D} is too weak. Feldvoss in \cite{F1} has replaced it with the condition $\mathrm{Im}(\delta_r) \subseteq \mathrm{Im}( 1 + \tau)$. We realized that the condition 
\begin{equation}\label{Hajjcond}
(\delta_r)^t(f\otimes f)=0 \text{ for all $f\in \fg^*$},
\end{equation}
which is weaker, also works.\\

\begin{mProposition}[Eq. \eqref{Hajjcond} $\Rightarrow r +\tau(r) \in (\fg\otimes \fg)^\fg$]
\textup{(}i\textup{)} We have: Eq. \eqref{Hajjcond} holds  if and only if 
\[
\sum_{1\leq i \leq n} \left ( f(b_i) f\circ \ad_{a_i}+ f(a_i) f \circ \ad_{b_i} \right )=0 \text{ for all $f\in \fg^*$.}
\]
\textup{(}ii\textup{)} If Eq. \eqref{Hajjcond} holds, then $r +\tau(r) \in (\fg\otimes \fg)^\fg$.

\end{mProposition}
\begin{proof}

Let us just prove part (ii). Observe that if Eq. (\ref{Hajjcond}) holds, then 
\[
(\delta_r)^t(f\otimes g+ g\otimes f)=0 \text{ for all $f,g\in \fg.$}
\] 
Now, for all $x\in \fg$, we have
\begin{align*}
0=& \langle f\otimes g + g \otimes f, \delta_r(x)\rangle = \langle f \otimes g, \delta_r(x) + \tau \circ \delta_r(x)\rangle.
\end{align*}
It remains to show that the condition 
\[
 \delta_r(x) + \tau \circ \delta_r(x)=0  \text{ for all $x\in \fg$}
\] 
is equivalent to $r +\tau(r) \in (\fg\otimes \fg)^\fg$. Indeed, 
\[
\delta_r(x) + \tau \circ \delta_r(x)=x\cdot r + \tau \circ x \cdot r= x \cdot (r + \tau(r)). \qed
\]
\noqed
\end{proof}


\begin{mTheorem}[Quasi-triangular Lie bialgebra structure,  see\cite{F1}]\label{Vald}
If $\fg$ is a~finite-dimen-sional non-abelian Lie algebra over an algebraically
closed field of characteristic $2$, then $\fg$ admits a~non-trivial quasi-triangular Lie
bialgebra structure.

\end{mTheorem}
\begin{proof} 
See \cite{F1}.
\end{proof}

\begin{mRemark} {\rm (i) If the Lie algebra is $\mathbb{Z}\slash 2$-graded and the quasi-triangular structure given by $r$ satisfies $r\in  (\fg \otimes \fg)_\ev$ then we say that $\fg$ admits {\it a~non-trivial quasi-triangular $\mathbb{Z}\slash 2$-graded Lie bialgebra structure.} 

(ii) If $\fg$ is a~Lie superalgebra in characteristic 2 where $\fg_\ev$ is a~non-abelian Lie algebra, then it is always possible to find a~non-trivial quasi-triangular $r$-matrix on $F(\fg)$. Indeed, just apply Theorem \ref{Vald} to $\fg_\ev$. 

(iii) In the case of Lie superalgebras over a~field of characteristic zero, and under certain conditions, several $r$-matrices have been investigated in \cite{Ka1, Ka2}.}\\
\end{mRemark}

Let $\fg=\fg_\ev\oplus \fg_\od$ be a~$\mathbb{Z}\slash 2$-graded vector space over $\mathbb K$. Let $c_\fg:\fg \rightarrow \fg\otimes\fg$ be an even linear map and $s_{\fg^*}:\fg_\od^* \rightarrow \fg_\ev^*$ be a map.  We say that $(\fg, c_\fg, s_{\fg^*})$ is a~Lie {\it supercoalgebra} if:

(i) $(\fg, c_\fg)$ is a~$\mathbb{Z}\slash 2$-graded Lie coalgebra; 

(ii) The map $s_{\fg^*}$ is a squaring on $\fg^*$ such that $(\fg^*, [\cdot, \cdot]_{\fg^*}, s_{\fg^*})$, where $ [\cdot, \cdot]_{\fg^*}$ is given as in \eqref{brag}, is a~Lie superalgebra.\\

Let $(\fg, [\cdot, \cdot]_\fg, s_\fg)$ be a~Lie superalgebra in characteristic $2$.  Let $c_\fg:\fg \rightarrow \fg\otimes\fg$ be an even linear map such that $(\fg, c_\fg, s_{\fg^*})$ is a~supercoalgebra in characteristic $2$. We say that $(\fg, [\cdot, \cdot]_\fg, c_\fg, s_\fg, s_{\fg^*})$ is a~Lie {\it superbialgebra} if:

(i) $c_\fg$ is a~1-cocyle on $\fg$ with values in $\fg\otimes \fg$; namely, it satisfies Eqs. (\ref{Cond1}) and (\ref{Cond2}).

(ii) $s_{\fg^*}$ satisfies the condition
\begin{equation}
\label{newsgstarbi}
s_{\fg^*}(f)\circ [\cdot, \cdot]_\fg= (f \otimes f)\circ (\id_\fg\otimes [\cdot, \cdot]_\fg)\circ (c_{\fg}\otimes \id_\fg+\id_\fg\otimes c_{\fg})\text{ for all $f\in \fg_1^*$}.
\end{equation}

As in the classical case, we will show that with our definitions the notion of a Manin triple is equivalent to the notion of superbialgebra. We have already shown in Proposition \ref{equicg} that if one  of the maps $c_\fg$ and $c_{\fg^*}$ satisfies the 1-cocycle condition (\ref{Cond1}),  then the other one must automatically satisfy this condition. It remains to show that the condition (\ref{newsgstarbi}) is equivalent to the 1-cocycle condition (\ref{Cond2}) for the map $c_{\fg^*}$. We have \\
\begin{mProposition}[Eq. (\ref{newsgstarbi}) is cohomological]
The squaring $s_{\fg^*}$ satisfies Eq. \eqref{newsgstarbi} if and and only if the map $c_{\fg^*}$ satisfies Eq. \eqref{Cond2} of the $1$-cocycle condition.

\end{mProposition}
\begin{proof} 

For all $f\in \fg_\od^*$, let us write $c_{\fg^*}(f)= \sum_{(f)} f_{(1)}\otimes f_{(2)}$. Now, for all $x,y\in \fg$ we have:
\begin{align*}
&\langle f\cdot c_{\fg^*}(f), x\otimes y \rangle=\langle \sum_{(f)} f\cdot f_{(1)}\otimes f_{(2)} + f_{(1)}\otimes f \cdot f_{(2)}, x \otimes y \rangle\\
&=\sum_{(f)} f\cdot f_{(1)}(x) f_{(2)}(y) + f_{(1)}(x) (f \cdot f_{(2)}) (y)\\
&=\sum_{(f)} (f\otimes f_{(1)})(c_\fg(x)) f_{(2)}(y) + f_{(1)}(x) (f \otimes f_{(2)}) (c_\fg(y))\\
&=\sum_{(f)} \langle f\otimes f_{(1)} \otimes f_{(2)}, c_\fg(x)\otimes y + x \otimes c_\fg(y) \rangle\\
&=\langle (f\otimes c_{\fg^*}(f)) \circ (c_\fg \otimes \id_\fg+ \id_\fg\otimes c_\fg), x\otimes y \rangle \\
& = \langle (f\otimes f )\circ (\id_\fg \otimes [\cdot, \cdot]_\fg)\circ (c_\fg \otimes \id_\fg+ \id_\fg \otimes c_\fg), x\otimes y\rangle.
\end{align*}
On the other hand
\[
\langle c_{\fg^*}(s_{\fg^*}(f)), x\otimes y \rangle =s_{\fg^*}(f)([x, y]_\fg) =\langle s_{\fg^*}(f) \circ [\cdot, \cdot]_\fg, x \otimes y\rangle. \qed
\]
\noqed\end{proof}

Let $\fg$ be a~non-abelian Lie superalgebra in characteristic $2$. Denote by $F(\fg)$ the Lie algebra structure on the underlying vector space obtained by forgetting the super structure. The Lie algebra $F(\fg)$ is $\mathbb{Z}\slash 2$-graded,  the grading being inherited from the super structure. Theorem \ref{Vald} guarantees the existence of a~non-trivial quasi-triangle Lie bialgebra structure on $F(\fg)$. Let us suppose further that this structure is $\mathbb{Z}\slash 2 \mathbb{Z}$-graded. Denote the corresponding $r-$matrix by $r$ and let us write 
 \[
r=\sum_{1 \leq i \leq n} a_i\otimes b_i , \quad \text{where \quad $\mathrm{gr}(a_i)=\mathrm{gr}(b_i)$ \quad (\text{namely, $r\in (\fg\otimes \fg)_\ev$).}}
\] 

Let us define the map $c_\fg:\fg\rightarrow \fg\otimes \fg$ by setting 
\begin{equation}\label{coboundF}
c_\fg(x):= \sum_{1 \leq i \leq n} \left (  [x,a_i]_\fg\otimes b_i + a_i\otimes [x,b_i]_\fg \right ) \quad  \text{for all  $x\in \fg$}.
\end{equation} 

\begin{mLemma}[$c_\fg$ is a 1-cocycle]
For all $x\in \fg_\od$, we have 
\[
x\cdot c_\fg(x)= c_\fg(s_\fg(x)).
\]

\end{mLemma}
\begin{proof} 

For all $x\in \fg_\od$, we have
\begin{align*}
& x\cdot c_\fg(x)= x\cdot \left (\sum_{1 \leq i \leq n} \left (  [x,a_i]\otimes b_i + a_i \otimes [x,b_i] \right ) \right )\\
&= \sum_{1 \leq i \leq n}  \left ( [x,[x,a_i]] \otimes b_i +[x, a_i] \otimes [x,b_i] \right )+  \sum_{1 \leq i \leq n} \left ( [x, a_i] \otimes [x,b_i]+ a_i\otimes [x,[x,b_i]] \right ) \\
&= \sum_{1 \leq i \leq n}  \left ( [s_\fg(x),a_i] \otimes b_i + a_i\otimes [s_\fg(x),b_i] \right ) = c_\fg(s_\fg(x)).\qed
\end{align*}
\noqed
\end{proof}

This lemma shows that the map $c_\fg$ defined by Eq. (\ref{coboundF}) is not only a~1-cocycle on the Lie algebra $F(\fg)$ but also a~1-cocyle on the Lie superalgebra $\fg$ as well. Therefore, we can write $c_\fg={\mathfrak d} r$, where $\mathfrak d$ is the coboundary operator,  see Eq. (\ref{diffsuper}).\\


\begin{mProposition}[A squaring on $\fg^*$ by means of an $r$-matrix] Let $\fg$ be a~Lie superalgebra in characteristic 2 with an $r$-matrix. There exists a~squaring on $\fg^*$ given as follows:
\begin{equation}\label{squgstar}
s_{\fg^*}: \fg^*_\od \rightarrow \fg^*_\ev \quad \quad  f\mapsto s_{\fg^*}(f):=\sum_{1 \leq i \leq n} f([-,a_i]_\fg) f(b_i).
\end{equation}

\end{mProposition}
\begin{proof} 

First, for all $\lambda\in \Kee$ and $f\in \fg^*_\od$ we have
\[
s_{\fg^*}(\lambda f)=\sum_{1 \leq i \leq n} \lambda f([-,a_i]_\fg) \lambda f(b_i)= \lambda^2 \sum_{1 \leq i \leq n} f([-,a_i]_\fg) f(b_i)=\lambda^2 s_{\fg^*}(f).
\]
Second, the map (for all $f,g\in \fg^*_\od$)
\[
\begin{array}{lcl}
(f,g) & \mapsto &  s_{\fg^*}(f+g)+s_{\fg^*}(f)+s_{\fg^*}(g)\\[1mm]
& &\displaystyle =\sum_{1 \leq i \leq n} \left (  f([-,a_i]_\fg) g(b_i)+ g([-,a_i]_\fg) f(b_i) \right ).
\end{array}
\]
is obviously bilinear.
\end{proof}
The existence of the squaring on $\fg^*$ does not endow $\fg^*$ with a~Lie superalgebra structure in characteristic $2$. To do that we will have to show that the Jacobi identity involving the squaring is satisfied. \\

\begin{mTheorem}[Condition on the squaring to satisfy the Jacobi identity]\label{sqgdual} Let $\fg$ be a~non-trivial Lie superalgebra in characteristic 2 such that $F(\fg)$ admits a~quasi-triangular Lie bialgebra structure defined by an even  tensor $r=\sum_{1\leq i \leq n}a_i \otimes b_i\in \fg\otimes \fg$. There exists a~Lie superalgebra structure on the dual space $\fg^*$, where the squaring is given by Eq. \eqref{squgstar}, if and only if the following condition is satisfied \textup{(}for all $x\in \fg$\textup{)}:
\begin{align}
\nonumber &\sum_{i,j} f(b_j) \left (f([[x,a_i]_\fg, a_j]_\fg)  h(b_i)+f([a_i, a_j]_\fg) h([x,b_i]_\fg) \right )\\[2mm]
\label {JIsqdual}&= \sum_{i,j} f([x,a_i]_\fg) ( f([b_i, a_j]_\fg h(b_j)+ f(a_j) h([b_i, b_j]_\fg)+ f(a_i) f([[x,b_i]_\fg, a_j]_\fg)h(b_j) \\[1mm]
\nonumber &+\sum_{i,j} f(a_i) f(a_j) h([[x,b_i]_\fg, b_j]_\fg ).
\end{align}

\end{mTheorem}
\begin{proof} 

We have already a~Lie algebra structure on $F(\fg^*)$ with the bracket given by 
\[
\begin{array}{lcl}
\langle [f,h]_{\fg^*}, x\rangle & = & \langle f\otimes h, c_\fg(x)\rangle\quad  \text{ for all $x\in \fg$ and $f,h\in \fg^*$}.
\end{array}
\]
The Jacobi identity involving the squaring reads
\begin{align*}
&\langle [s_{\fg^*}(f),h]_{\fg^*},x \rangle= \langle s_{\fg^*}(f) \otimes h, c_\fg(x)\rangle= \langle s_{\fg^*}(f) \otimes h, \sum_{1\leq i \leq n} \left (  [x,a_i]_\fg \otimes b_i + a_i \otimes [x,b_i]_\fg \right ) \rangle\\
&= \sum_{1\leq i \leq n} \left (   s_{\fg^*}(f)( [x,a_i]_\fg) h(b_i) + s_{\fg^*}(f)(a_i) h([x,b_i]_\fg) \right ) \\
&=   \sum_{1\leq i,j \leq n} \left (  f([[x,a_i]_\fg ,a_j]_\fg) f(b_j) h(b_i) + f( [a_i,a_j]_\fg) h([x,b_i]_\fg) \right ).
\end{align*}
On the other hand, 
\begin{align*}
&\langle [f,[f,h]_{\fg^*}]_{\fg^*},x \rangle=   \langle f \otimes [f,h]_{\fg^*}, c_\fg(x)\rangle  =    \langle f \otimes [f,h]_{\fg^*},  \sum_{1\leq i \leq n} \left (  [x,a_i]_\fg \otimes b_i + a_i \otimes [x,b_i]_\fg \right ) \rangle  \\
&= \sum_{1\leq i \leq n} \left (  f([x,a_i]_\fg) [f,h]_{\fg^*}(b_i) +f(a_i)  [f,h]_{\fg^*}([x,b_i]_\fg) \right )\\
&= \sum_{1\leq i \leq n} \left (  f([x,a_i]_\fg) \langle f\otimes h, c_\fg(b_i)\rangle +f(a_i)  \langle f \otimes h, c_\fg([x,b_i]_\fg)\rangle \right )\\
&= \sum_{1\leq i,j \leq n}  f([x,a_i]_\fg) (  f([b_i,a_j]_\fg)  h(b_j) + f(a_j) h([b_i,b_j]_\fg) )\\
& + \sum_{1\leq i,j \leq n} f(a_i)   \left (  f([ [x,b_i]_\fg ,a_j]_\fg) h(b_j) + f(a_j) h([ [x,b_i]_\fg ,b_j]_\fg) \right ).\qed
\end{align*}
\noqed\end{proof}
We have just proved that if the condition (\ref{JIsqdual}) is satisfied, then we have a structure of Lie supercogebra on $(\fg, c_\fg, s_{\fg^*})$, where $s_{\fg^*}$ is defined by Eq. (\ref{squgstar}). 

As a result, we reach the following definition.  A classical $r$-matrix is called {\it admissible} if it is even and satisfies the  condition (\ref{JIsqdual}). 
\begin{mRemark}
{\rm The admissibility condition guarantees the existence of a squaring on $\fg^*$ satisfying the Jacobi identity. }\\
\end{mRemark}

Let us now study the property of the map $c_{\fg^*}$ as given in Eq. (\ref{bragstar}). We will show that it satisfies the 1-cocycle condition involving the squaring.  First, we need the following lemma.\\

\begin{mLemma}[Technical] \label{Eqh}  For all $x,y\in \fg$ and for all $f\in \fg^*$, we have
\[
f\cdot c_{\fg^*}(f)(x\otimes y)=(f\otimes y\cdot f)(c_\fg(x))+(f\otimes x\cdot f)(c_\fg(y)).
\]

\end{mLemma}
\begin{proof}

Indeed, 
\begin{align*}
& f\cdot c_{\fg^*}(f)(x\otimes y) =  c_{\fg^*}(f)(f\cdot x \otimes y)+c_{\fg^*}(f)(x\otimes f\cdot y)= f([f\cdot x, y]_\fg)+f([x,f\cdot y]_\fg)\\[2mm]
&=f([x\circ \ad_f,y]_\fg)+f ([x,y\circ \ad_f]_\fg) =(y\cdot f)(x\circ \ad_f)+(x\cdot f)(y\circ \ad_f)\\[2mm]
&=[f, y\cdot f]_{\fg^*}(x)+[f,x\cdot f]_{\fg^*}(y) =(f\otimes y\cdot f)(c_\fg(x))+(f\otimes x\cdot f)(c_\fg(y)). \qed
\end{align*}
\noqed\end{proof}

\begin{mProposition}[$c_{\fg^*}$ is a 1-cocycle]
 We have $c_{\fg^*}\in Z^1(\fg^*, \fg^*\otimes \fg^*)$. 

\end{mProposition}
\begin{proof}

Let us check that condition Eq. (\ref{Cond2}) is satisfied. Using Lemma \ref{Eqh} we get (for all $x,y\in \fg$ and for all $f\in \fg^*_\od$):
\[
\begin{array}{lcl}
f\cdot c_{\fg^*}(f)(x\otimes y)&=&(f\otimes y\cdot f)(\sum_{1\leq i \leq n} [x,a_i]_\fg\otimes b_i+a_i\otimes [x,b_i]_\fg)+(x\leftrightarrow y)\\[2mm]
&=&\sum_{1\leq i \leq n} f([x,a_i]_\fg) f([y,b_i])+f(a_i)f([y,[x,b_i]_\fg]_\fg) +(x\leftrightarrow y) \\[2mm]
&=&\sum_{1\leq i \leq n} f(a_i) f([[x,y]_\fg,b_i]_\fg)= c_{\fg^*} (s_{\fg^*}(f))(x\otimes y). \qed
\end{array}
\]
\noqed\end{proof}
\begin{mCorollary}
Any admissible $r$-matrix on $\fg$ yields a Lie bisuperalgebra structure on $\fg$.  \\
\end{mCorollary}

Consider $r \in \fg \otimes \fg$ given by $r=\sum_{1\leq i \leq n} a_i\otimes b_i$, and define a~linear map $R:\fg^*\rightarrow \fg$ by setting
\begin{equation}
\label{defR}
R(f):=(f\otimes 1)(r)=\sum_{1\leq i \leq n} f(a_i) b_i \quad \text{ for all $f\in \fg^*$.}
\end{equation}
In the sequel,  we will write the properties of the $r$-matrix $r$ in terms of the linear map $R$ (called and $R$-matrix).  

\begin{mProposition}[Technical] \label{IJR}
Let $\fg$ be a~Lie superalgebra in characteristic $2$, and let $R$ be the map \eqref{defR}  defined by means of $r\in \fg\otimes \fg$. Then, $r$ is an even and symmetric classical $r$-matrix if and only if $R$ is even and satisfies  the following two conditions:\\

\textup{(}i\textup{)} $\langle f, R(h)\rangle =   \langle h, R(f)\rangle \quad \text{for all $f,h\in \fg^*$};$ \\

\textup{(}ii\textup{)} $\langle f, [R(h), R(l)]_\fg\rangle + \langle l, [R(f), R(h)]_\fg\rangle +  \langle h, [R(l), R(f)]_\fg\rangle=0 \quad \text{for all $f,h,l\in \fg^*$.}$

\end{mProposition}
\begin{proof} 

Obviously, $r$ is even if and only if $R$ is even. Also $r$ is symmetric if and only if Part (i) is satisfied. Now, a~direct computation shows that 
\[
\begin{array}{l}
\langle f, [R(h), R(l)]_\fg\rangle + \langle l, [R(f), R(h)]_\fg\rangle +  \langle h, [R(l), R(f)]_\fg\rangle  \\[1mm]
= \langle f\otimes h \otimes l, [r_{12},r_{13}]+[r_{12}, r_{23}]+[r_{13}, r_{23}]\rangle.\qed
\end{array}
\]  \noqed \end{proof}

An even linear map $R:\fg^*\rightarrow \fg$ that satisfies (i) and (ii) of Theorem \ref{IJR} is called a classical $R$-matrix. 

\begin{mTheorem}[A bracket and a squaring on $\fg^*$]
Let $\fg$ be a~Lie superalgebra in characteristic $2$.  Let $F(\fg)$ admit a~triangular Lie bi-algebra structure such that the co-multiplication on $F(\fg)$ is the $1$-coboundary of $r$, where $r$ is an even symmetric $r$-matrix. Then the multiplication on $\fg^*$ defined by Eq. \eqref{bragstar} is explicitly determined by the formula 
\[
[f,h]_{\fg^*}=h \circ \ad_{R(f)}+ f \circ \ad_{R(h)}\quad \text{for all $f,h\in \fg^*$.}
\]
and the squaring given by Eq. \eqref{squgstar} is explicitly given by the formula 
\begin{equation}\label{squgstar2}
s_{\fg^*}(f)=f \circ \ad_{R(f)}.
\end{equation}

\end{mTheorem}
\begin{proof}

Let us first deal with the bracket.  For all $x\in \fg$, the LHS reads
\begin{align*}
&\langle [f,h]_{\fg^*},x\rangle=\langle f\otimes h,c_\fg(x)\rangle  =\langle f\otimes h, \sum_{1 \leq i \leq n} ([x,a_i]\otimes b_i + a_i \otimes [x,b_i] ) \rangle \\
&= \sum_{1 \leq i \leq n} f([x,a_i]) h(b_i) +f(a_i) h([x,b_i]).
\end{align*} 
The RHS reads
\begin{align*}
&\langle h \circ \ad_{R(f)}+ f \circ \ad_{R(h)}\,,x\rangle=h \circ [R(f),x]_\fg+ f \circ [R(h),x]_\fg= \sum_{1 \leq i \leq n} f([x,a_i]) h(b_i) +f(a_i) h([x,b_i]).
\end{align*} 
When it comes to the squaring, and using the fact that $r$ is symmetric, we get:
\begin{align*}
& \langle f \circ \ad_{R(f)}, x \rangle=f([R(f),x]_\fg)= f ([\sum_{1 \leq i \leq n} f(a_i) b_i,x]_\fg)\\
&= \sum_{1 \leq i \leq n} f(a_i) f ([b_i,x]_\fg) =\sum_{1 \leq i \leq n} f(b_i) f ([a_i,x]_\fg) 
=\langle s_{\fg^*}(f), x \rangle. \qed
\end{align*}
\noqed\end{proof}

\begin{mProposition}[A Lie superalgebra structure on $\fg^*$]\label{sqgdual2}
Let $\fg$ be a~Lie superalgebra in characteristic $2$. Let $F(\fg)$ admit a~triangular Lie bi-algebra structure such that the co-multiplication on $F(\fg)$ is the $1$-coboundary of $r$ \textup{(}where $r$ is an even symmetric $r$-matrix\textup{)}. Then,  there exists a~Lie superalgebra structure on the dual space $\fg^*$, where the squaring is given by \eqref{squgstar2}, if and only if the following condition is satisfied \textup{(}for all $x\in \fg$, for all $h\in \fg^*$ and for all $f\in \fg^*_\od$\textup{)}:
\begin{equation} \label{eq26b}
\begin{array}{l}
 h \circ \ad_{s_\fg(R(f))}+ f \circ \ad_{[R(h), R(f)]_\fg}+ f \circ \ad_{R(h \circ \ad_{R(f)})+R( f \circ \ad_{R(h)})}+h\circ \ad_{R(f\circ \ad_{R(f)})}=0.
\end{array}
\end{equation}

\end{mProposition}
\begin{proof}

For the proof of the Jacobi identity on $\fg^*_\ev$ and the fact that $\fg_\od$ is a~$\fg_\ev$-module, we refer to \cite[Prop. 3.1]{BaB}. Let us deal with the squaring
\begin{align*}
& [s_{\fg^*}(f),h] _{\fg^*}= h\circ \ad_{R(s_{\fg^*}(f))}+ s_{\fg^*}(f) \circ \ad_{R(h)}=h\circ \ad_{R(f\circ \ad_{R(f)})}+ f\circ \ad_{R(f)}\circ \ad_{R(h)}.
\end{align*}
On the other hand,
\begin{align*}
&[f,[f,h]_{\fg*}]_{\fg*}=[f,h \circ \ad_{R(f)}+ f \circ \ad_{R(h)}]_{\fg*}\\[1mm]
&= (h \circ \ad_{R(f)}+ f \circ \ad_{R(h)}) \circ \ad_{R(f)}+ f \circ \ad_{R(h \circ \ad_{R(f)})+R( f \circ \ad_{R(h)})}\\[1mm]
& = h \circ \ad_{s_\fg(R(f))}+ f \circ \ad_{R(h)} \circ \ad_{R(f)}+ f \circ \ad_{R(h \circ \ad_{R(f)})+R( f \circ \ad_{R(h)})}
 \qed
\end{align*}
\noqed\end{proof}
\begin{mRemark}{\rm 
In fact, Eq. (\ref{JIsqdual}) is equivalent to Eq. (\ref{eq26b}).}\\
\end{mRemark}

Consider the space
\[
\mathrm{Im}(R):=\{ R(f)\; | \; f \in \fg^*\}, 
\]
and consider the following bilinear form
\begin{equation}\label{ome}
\omega:\mathrm{Im}(R) \times \mathrm{Im}(R)\rightarrow \mathbb{K} \quad (R(f), R(g))\mapsto \langle g, R(f)\rangle.
\end{equation}
The form $\omega$ is well-defined. Indeed, if $R(f)=R(f')$ and $R(h)=R(h')$, then
\[
\langle f, R(h)\rangle = \langle R(f), h \rangle =\langle R(f'), h\rangle= \langle f', R(h)\rangle= \langle f', R(h')\rangle.
\]
\begin{mProposition}[A closed non-degenerate $\od$-antisymmetric form on $\mathrm{Im}(R)$] The space 
$\mathrm{Im}(R)$ is a $\mathbb Z/2$-graded Lie algebra in characteristic $2$ equipped with the form \eqref{ome} that is closed,  non-degenerate and $\od$-antisymmetric.

\end{mProposition}
\begin{proof}

We should prove first that $\mathrm{Im}(R)$ is closed under the bracket,  and hence $\mathrm{Im}(R)$ is a~Lie sub-algebra of $F(\fg)$. Indeed, let $R(f), R(h)\in \mathrm{Ker}(R)$. For all $l\in \mathrm{Ker}(R)$,  using Part (ii) of Prop. \ref{IJR},  we see that 
\[
\langle l, [R(f), R(h)]_\fg\rangle= \langle h, [R(l), R(f)]_\fg\rangle + \langle f, [R(h), R(l)]_\fg\rangle=0.
\]
It follows that $[R(f), R(h)]_\fg \in (\mathrm{Ker}(R))^\perp$, where here the orthogonality should be understood with respect to the bilinear form $\langle-,-\rangle$. It follows that  $[R(f), R(h)]_\fg\in \mathrm{Im}(R)$.

It is easy to see that the map $\omega$ is bilinear.  Let us check the 2-cocycle condition. Indeed, using Part (ii) of Proposition \ref{IJR}, we get
\begin{align*}
&\omega(R(l),[R(f),R(h)]_\fg)+\omega(R(h),[R(l),R(f)]_\fg)+\omega(R(f),[R(h),R(l)]_\fg)\\
&=\langle l, [R(f),R(h)]_\fg) \rangle +\langle h, [R(l),R(f)]_\fg) \rangle+ \langle f, [R(h),R(l)]_\fg) \rangle=0.
\qed
\end{align*}
\noqed
\end{proof}

Suppose now that $\fg$ has an~even NIS ${\mathscr B}$. Consider the natural isomorphism
\[
\Phi: \fg \rightarrow \fg^*  \quad x\mapsto \Phi(x):={\mathscr B}(x,\cdot).
\]
Set a~new bracket and a~new squaring on $\fg$ as follows:
\[
\begin{array}{lcll}
\widetilde{[\cdot, \cdot]} & : & \fg \times \fg \rightarrow \fg &  \widetilde{[x, y]}=\Phi^{-1}([\Phi(x), \Phi(y)]_{\fg^*}) \text{ for all $x,y\in \fg$}\\[1mm]
\widetilde{s_\fg(\cdot)} & : & \fg_\od \rightarrow \fg_\ev &  \widetilde{s_\fg(x)}=\Phi^{-1}(s_{\fg^*} (\Phi(x))) \text{ for all $x\in \fg_\od$}.
\end{array} 
\]

\begin{mTheorem}[Deforming the Lie structure] Let $(\fg, [\cdot, \cdot]_\fg, s_\fg )$ be a~Lie superalgebra in characteristic $2$ with an~even NIS. Suppose that $\fg$ has an invertible,   symmetric and  admissible $r$-matrix. 
If we denote $U:=R\circ \Phi:\fg \rightarrow \fg$,  then 
\[
\begin{array}{lcl}
\widetilde{[x,y]} & = & [U(x), y]_\fg+[x, U(y)]_\fg\quad \text{for all $x,y \in \fg$,}\\[2mm]
\widetilde{s}_\fg(x) & = & [U(x), x]_\fg \quad \text{for all $x\in \fg_\od$}.
\end{array}
\]
Moreover, $(\fg, \widetilde{[\cdot, \cdot]}, \widetilde{s}_\fg)$ is a~Lie superalgebra.

\end{mTheorem}
\begin{proof}
As the $r$-matrix is admissible,  this guarantees the existence of a Lie superalgebra structure on $\fg^*$.  We will only be dealing with the squaring. Indeed, for all $x\in \fg_\od$ we have
\begin{align*}
&\widetilde{s_\fg(x)}= \Phi^{-1}(s_{\fg^*} (\Phi(x))) = \Phi^{-1}( \Phi(x) \circ \ad_{R(\Phi(x))})= \Phi^{-1}( {\mathscr B}(x, [R(\Phi(x)), \cdot]_\fg))\\
&= \Phi^{-1}( {\mathscr B}([x, R(\Phi(x))]_\fg, \cdot))=[x, R(\Phi(x))]_\fg=[x, U(x)]_\fg.
\end{align*}
The JI pertaining to the squaring reads: (for all $x\in\fg_\od$ and for all $y\in \fg$)
\begin{align*}
&\widetilde{[\widetilde{s}_\fg(x), y]}= \Phi^{-1}([s_{\fg^*}(\Phi(x)), \Phi(y)]_{\fg^*})=  \Phi^{-1}([\Phi(x),[\Phi(x), \Phi(y)]_{\fg^*}]_{\fg^*}) =   \Phi^{-1}([\Phi(x), \Phi (\widetilde{[x, y]})]_{\fg^*})\\
& = \widetilde{[x, \widetilde{[x, y]}]}.\qed
\end{align*}
\noqed
\end{proof}
\begin{mProposition}[Properties of  the map $U$]\label{propUZ}  Let $(\fg, [\cdot, \cdot]_\fg, s_\fg )$ be a~Lie superalgebra in characteristic $2$ with an~even NIS.  Let $r\in \fg\otimes \fg$ that is symmetric and even, and let us put $U:=R\circ \Phi:\fg \rightarrow \fg$. 
\begin{itemize}
\item[(i)] Part \textup{(}i\textup{)} of Prop. \ref{IJR} holds if and only if $U \text{ is ${\mathscr B}$-symmetric}$.
\item[(ii)] Part \textup{(}ii\textup{)} of Prop. \ref{IJR} holds  if and only if 
\[
U(\widetilde{[x,y]})=[U(x), U(y)]_\fg \quad  \text{ for all $x,y \in \fg$}.
\]
\item[(iii)] Eq.  \eqref{eq26b} holds if and only if \textup{(}for all $x\in \fg_\od$ and $y\in \fg$\textup{)}:
\[
\ad_y(U(\tilde s_\fg(x)) + s_\fg(U(x)))=\ad_x (U(\widetilde{[x,y]})+[U(x),U(y)]_\fg).
\]
\end{itemize}
\end{mProposition}
\begin{proof}
We will prove only Part (i) and (iii). Let us start with Part (i). For all $f,g \in \fg^*$,  let us write $f=\Phi(a)$ and $g=\Phi(b)$, for some $a,b\in \fa$.  We get
\[
\langle g, R(f) \rangle= \langle \Phi(b),  R\circ \Phi(a) \rangle= \langle \Phi(b), U(a) \rangle= {\mathscr B}(b,U(a)).
\]
On the other hand, 
\[
\langle f, R(g) \rangle= \langle  \Phi(a),  R\circ \Phi(b) \rangle= \langle \Phi(a), U(b) \rangle= {\mathscr B}(a, U(b)).
\]
Let us prove Part (iii). For all $f\in \fg_\od$ and $g \in \fg$, let us write $f=\phi(a)$ and $g=\phi(b)$ for $a\in \fg_\od$ and $b\in \fg$. We have 
\begin{align*}
& h \circ \ad_{s_\fg(R(f))}+ f \circ \ad_{[R(h), R(f)]_\fg}+ f \circ \ad_{R(h \circ \ad_{R(f)})+R( f \circ \ad_{R(h)})}+h\circ \ad_{R(f\circ \ad_{R(f)})}=\\[1mm]
&=\Phi(b)\circ \ad_{s_\fg(U(a))}+ \Phi(a) \circ ( \ad_{[U(b), U(a)]_\fg}+  \circ \ad_{R(\Phi(b) \circ \ad_{U(a)}))+R( \Phi(a) \circ \ad_{U(b)})}+ \Phi(b)\circ \ad_{R(\Phi(a)\circ \ad_{U(a)})}\\[1mm]
&=\Phi([b,s_\fg(U(a))]_\fg)+\Phi([a, [U(a),U(b)]_\fg]_\fg)+\Phi(a)\circ \ad_{U([b,U(a)]_\fg)+U([a,U(b)]_\fg) }+ \Phi(b)\circ \ad_{U([a,U(a)]_\fg)}\\[1mm]
&=\Phi \left ([b,s_\fg(U(a))]_\fg + [a, [U(a),U(b)]_\fg]_\fg + [a, U([b,U(a)]_\fg)+U([a,U(b)]_\fg )]_\fg + [b, U([a,U(a)]_\fg) ]_\fg \right )\\[1mm]
&=\Phi \left (\ad_b( s_\fg(U(a))+U(\tilde s_\fg(a))+[a, U(\widetilde{[a,b]})+[U(a),U(b)]_\fg]_\fg \right ),
\end{align*}
using the  fact that $\tilde s_\fg(a)=[a, U(a)]_\fg$.
\end{proof}
\begin{mCorollary}
Let $(\fg, [\cdot, \cdot]_\fg, s_\fg )$ be a~Lie superalgebra in characteristic $2$ with an~even NIS $\mathscr B$.  Let $r\in \fg\otimes \fg$ that is symmetric and even, and let us put $U:=R\circ \Phi:\fg \rightarrow \fg$.  Then $r$ is an $r$-matrix if and only if $U$ is $\mathscr B$-invariant and $U:(F(\fg), \widetilde{[\cdot, \cdot]} )\rightarrow (F(\fg), [\cdot, \cdot]_\fg ) $ is a homomorphism of Lie algebras.  Moreover, if $r$ is an $r$-matrix, then $r$ is admissible if and only if
\[
\ad_y(U(\tilde s_\fg(x)) + s_\fg(U(x)))=0\text{ for all $x\in \fg_\od$ and $y\in \fg$.}
\]
\end{mCorollary}
\begin{mCorollary}
Let $(\fg, [\cdot, \cdot]_\fg, s_\fg )$ be a Lie superalgebra in characteristic $2$ with an~even NIS such that $\fz(F(\fg))=\{0\}$.  Let $r\in \fg\otimes \fg$ that is symmetric, even and admissible. Let us put $U:=R\circ \Phi:\fg \rightarrow \fg$.  Then the map $U:(\fg, \widetilde{[\cdot, \cdot]}, \tilde s_\fg )\rightarrow (\fg,  [\cdot, \cdot]_\fg , s_\fg) $ is a homomorphism of Lie superalgebras.
\end{mCorollary}
\begin{mRemark}
We will see in the last section that the condition $\fz(F(\fg))=\{0\}$ is not necessary to guarantee that $U$ is a homomorphism of Lie superalgebras. 
\end{mRemark}
\begin{mProposition}[A converse of Proposition \ref{propUZ}]\label{propU}
Let $U\in \End(\fg)$ be an even map such that ${\mathscr B}(U(x),y)={\mathscr B}(x, U(y))$,  and  
\[
\begin{array}{lcl}
\widetilde{[x,y]}_\fg & := & [U(x),y]_\fg+ [x, U(y)]_\fg \quad \text{for all $x,y\in \fg$}\\[2mm]
\widetilde{s}_\fg(x) & := & [x,U(x)]_\fg\quad \text{for all $x\in \fg_\od$}
\end{array}
\] 
defines a~Lie superalgebra structure on $\fg$.  Then, 
\[
R:=U\circ \Phi^{-1}:\fg^* \rightarrow \fg  \text{ is a~associated with an admissible $r$-matrix.}
\]

\end{mProposition}
\begin{proof}

The part not involving the squaring is similar to that of \cite[Prop 3.10 and  Prop 3.11]{BaB}. The JI involving the squaring $\tilde s_\fg$ is equivalent to the fact that $R$ verifies Eq. (\ref{eq26b}). \end{proof}

As a result, we arrive at the following definition. A map $U\in \End(\fg)$ satisfying the condition of Prop. \ref{propU} is called an {\it  admissible Yang-Baxter} operator (associated with an admissible $r$-matrix $r$). 


\subsection{The matched pairs for Lie superalgebras in characteristic $2$}\label{mpairs}
The goal of this subsection is to look at Manin triples in a~more general setting by introducing the concept of matched pairs. As far as we know, this concept was introduced by Majid in \cite{Maj}.

Let $(\fg,[\cdot, \cdot]_\fg,s_\fg)$ and $(\fk, [\cdot, \cdot]_\fk, s_\fk)$ be two Lie superalgebras in characteristic $2$ -- in the case where  $\fk=\fg^*$,  the triple $(\fh:=\fg\oplus \fg^*, \fg, \fg^*)$ is a Manin triple.  Let us suppose they are both endowed with the representations
\[
\pi:\fg_\ev \rightarrow \mathrm{End}(\fk_\ev), \quad \lambda:\fg_\ev \rightarrow \mathrm{End}(\fk_\od), \quad \rho:\fk_\ev \rightarrow \mathrm{End}(\fg_\ev), \quad \mu:\fk_\ev \rightarrow \mathrm{End}(\fg_\od),
\]
together with the following maps
\[
\tilde \lambda: \fg_\ev \rightarrow \mathrm{Hom} (\fk_\od, \fg_\od), \quad \tilde \mu: \fk_\ev \rightarrow \mathrm{Hom} (\fg_\od, \fk_\od).
\]
Now suppose that these maps satisfy the following conditions: (for all $x_0,y_0\in\fg_\ev$, $x_1\in\fg_\od$, $a_0,b_0\in\fk_\ev$, $a_1\in\fk_\od$):
\begin{equation}
\label{hev}
\begin{array}{l}
\rho(a_0) ([x_0, y_0]_\fg)+ [\rho (a_0)(x_0), y_0]_\fg+[x_0, \rho (a_0)(y_0)]_\fg\\[1mm]
=\rho(\pi(y_0)(a_0))(x_0)+\rho(\pi(x_0)(a_0))(y_0),\\[2mm]
 \pi(x_0) ([a_0, b_0]_\fk)+ [\pi (x_0)(a_0), b_0]_\fk+[a_0, \pi (x_0)(b_0)]_\fk\\[1mm]
=\pi(\rho(b_0)(x_0))(a_0)+\pi(\rho(a_0)(x_0))(b_0), 
 \end{array}
 \end{equation}
 and
 \begin{equation}
 \begin{array}{lcl}
 \label{hod}
 \mu(b_0) ([x_0, x_1]_\fg)+[x_0, \mu (b_0)(x_1)]_\fg&=&[\rho(b_0)(x_0),x_1]_\fg+\mu(\pi(x_0)(b_0))(x_1) + \tilde \lambda (x_0) (\tilde \mu (b_0)(x_1)),\\[3mm]
\lambda(x_0) ([b_0, a_1]_\fk)+[b_0, \lambda (x_0)(a_1)]_\fk&=&[\pi(x_0)(b_0),a_1]_\fk+\lambda(\rho(b_0)(x_0))(a_1) +\tilde \mu (b_0) (\tilde \lambda (x_0)(a_1)),\\[3mm]
 \tilde \lambda (\rho(b_0)(x_0))(a_1)+\tilde \lambda (x_0)([b_0,a_1]_\fk)&=&\mu (b_0)(\tilde \lambda (x_0)(a_1)),\\[3mm]
  \tilde \mu (\pi(x_0)(b_0))(x_1)+\tilde \mu (b_0)([x_0,x_1]_\fg)&=&\lambda (x_0)(\tilde \mu (b_0)(x_1)),\\[3mm]
 \tilde  \lambda([x_0,y_0]_\fg) (a_1)+[x_0, \tilde \lambda (y_0)(a_1)]_\fg&=&[y_0, \tilde \lambda (x_0)(a_1)]_\fg+\tilde \lambda(x_0) \lambda(y_0)(a_1)+\tilde \lambda (y_0)  \lambda (x_0)(a_1),\\[3mm]
  \tilde  \mu([a_0,b_0]_\fk) (x_1)+[a_0, \tilde \mu (b_0)(x_1)]_\fk&=&[b_0, \tilde \mu (a_0)(x_1)]_\fk+\tilde \mu(a_0) \mu(b_0)(x_1) +\tilde \mu (b_0)  \mu (a_0)(x_1).
\end{array}
\end{equation}
Suppose further that there are two bilinear maps $r_\fg: \fg_\od \odot  \fk_\od \rightarrow \fg_\ev$ and  $r_\fk: \fg_\od\odot \fk_\od \rightarrow \fk_\ev$, satisfying the following conditions (for all $b\in \fk_\ev$ and for all $y\in \fg_\ev$):

\begin{equation}
\label{sq5}
\begin{array}{lcl}
 [r_\fg(x_1,a_1),y]_\fg+\rho (r_\fk(x_1,a_1))(y)&=&r_\fg(x_1, \lambda(y) (a_1))+[x_1, \tilde \lambda(y)(a_1)]_\fg+r_\fg(a_1,[x_1,y]_\fg),\\[3mm]
\rho(b)(r_\fg(x_1, a_1))+r_\fg(a_1, \mu(b)(x_1))&=&r_\fg(x_1, [a_1,b]_\fk),\\[3mm]
\rho(b)(s_\fg(x_1))+[x_1, \mu(b)(x_1)]_\fg&=& r_\fg(x_1, \tilde \mu (b)(x_1)),\\[3mm]
\rho(s_\fk(a_1))(y)&=& r_\fg(a_1, \tilde \lambda (y)(a_1)),
\end{array}
\end{equation}
and also
\begin{equation}
\label{sq6}
\begin{array}{lcl}
 [r_\fk(x_1,a_1),b]_\fk+\pi (r_\fg(x_1,a_1))(b)&=&r_\fk(a_1, \mu(b) (x_1))+[a_1, \tilde \mu (b)(x_1)]_\fk+r_\fk(x_1,[a_1,b]_\fk),\\[3mm]
\pi(y)(r_\fk(x_1, a_1))+r_\fk(x_1, \lambda(y)(a_1))&=&r_\fk(a_1, [x_1,y]_\fg),\\[3mm]
\pi(y)(s_\fk(a_1))+[a_1, \lambda(y)(a_1)]_\fk&=& r_\fk(a_1, \tilde \lambda (y)(a_1)),\\[3mm]
\pi(s_\fg(x_1))(b)&=& r_\fk(x_1, \tilde \mu (b)(x_1)).
\end{array}
\end{equation}
Moreover, (for all $b\in \fk_\od$ and for all $y\in \fg_\od$)
\begin{equation}
\label{sq7}
\begin{array}{lcl}
\tilde \mu (r_\fk(x_1, a_1))(y)&=&\lambda ([x_1, y]_\fg)(a_1)+\tilde \mu (r_\fk(a_1, y))(x_1),\\[3mm]
\tilde \mu (s_\fk(a_1))(y)&=& [r_\fk(a_1, y),a_1]_\fk+\lambda (r_\fg(a_1,y))(a_1),\\[3mm]
\lambda (s_\fg(x_1))(b)&=&\tilde \mu (r_\fk(x_1, b))(x_1),\\[3mm]
\lambda(r_\fg(x_1, a_1))(b)&=&[r_\fk(x_1, a_1),b]_\fk+\lambda (r_\fg(x_1, b))(a_1) +\tilde \mu ([a,b]_\fk)(x_1)+[r_\fk(x_1, b),a_1]_\fk,
\end{array}
\end{equation}
also
\begin{equation}
\label{sq8}
\begin{array}{lcl}
\tilde \lambda (r_\fg(x_1, a_1))(b)&=&\mu ([a, b]_\fk)(x_1)+\tilde \lambda (r_\fg(x_1, b))(a),\\[3mm]
\tilde \lambda (s_\fg(x_1))(b)&=& [r_\fg(x_1, b),x_1]_\fg+\mu (r_\fk(x_1,b))(x),\\[3mm]
\mu (s_\fk(a_1))(y)&=&\tilde \lambda (r_\fg(a_1, y))(a_1),\\[3mm]
\mu(r_\fk(x_1, a_1))(y)&=&[r_\fg(x_1, a_1),y]_\fg+\mu (r_\fk(a_1, y))(x_1) +\tilde \lambda (a_1)([x_1,y]_\fg)+[r_\fg(a_1, y),x_1]_\fg.
\end{array}
\end{equation}

\begin{mTheorem}[Matched pairs] \label{mainthm} Let $(\fg,[\cdot, \cdot]_\fg,s_\fg)$ and $(\fk, [\cdot, \cdot]_\fk, s_\fk)$ be two Lie superalgebras in characteristic $2$, endowed with 
\begin{itemize}
\item Representation maps
\[
\pi:\fg_\ev \rightarrow \mathrm{End}(\fk_\ev), \quad \lambda:\fg_\ev \rightarrow \mathrm{End}(\fk_\od), \quad \rho:\fk_\ev \rightarrow \mathrm{End}(\fg_\ev), \quad \mu:\fk_\ev \rightarrow \mathrm{End}(\fg_\od);
\]
\item Linear maps $\tilde \lambda: \fg_\ev \rightarrow \mathrm{Hom} (\fk_\od, \fg_\od)$ and $\tilde \mu: \fk_\ev \rightarrow \mathrm{Hom} (\fg_\od, \fk_\od)$ satisfying Eqns \eqref{hev}, \eqref{hod};
\item Bilinear maps $r_\fg: \fg_\od\odot \fk_\od \rightarrow \fg_\ev$ and  $r_\fk: \fg_\od\odot \fk_\od \rightarrow \fk_\ev$, satifying  Eqns \eqref{sq5}, \eqref{sq6}, \eqref{sq7}, \eqref{sq8}.
\end{itemize}
Then, there exists a~Lie superalgebra structure on $\fh:=\fg\oplus \fk$ given as follows.
\begin{itemize}
\item \underline{ The Lie algebra structure on $\fh_\ev$:} The Lie bracket is given by:
\[
[x_0+a_0, y_0+b_0]_\fh:=[x_0,y_0]_\fg+ \rho(b_0)(x_0)+ \rho(a_0)(y_0)+[a_0,b_0]_\fk+ \pi(x_0)(b_0)+ \pi(y_0)(a_0),
\]
for all $x_0, y_0 \in \fg_\ev$ and $a_0,b_0 \in \fk_\ev$.
\item \underline{ The $\fh_\ev$-module structure on $\fh_\od$:} The action is given by:
\[
[x_0+a_0, x_1+a_1]:=[x_0,x_1]_\fg+  \mu(a_0)(x_1)+ \tilde \lambda(x_0)(a_1)+[a_0,a_1]_\fk+\lambda(x_0)(a_1)+\tilde \mu(a_0)(x_1),
\]
for all $x_0\in \fg_\ev, x_1\in \fg_\od$ and $a_0 \in \fk_\ev, b_1 \in \fk_\od$.
\item \underline{ The squaring on $\fh_\od$} \textup{(}where we have put $r:=r_\fg+r_\fk$\textup{)}:
\[
S_\fh(x_1+a_1):=s_\fg(x_1)+s_\fk(a_1)+r(x_1,a_1),
\]
for all $x_1\in \fg_\od$ and $a_1 \in \fk_\od$.
\end{itemize}

\end{mTheorem}
\begin{proof}

The Jacobi identity for $\fh_\ev$ is equivalent to the fact that $\pi$ and $\rho$ are representations satisfying  conditions (\ref{hev}). We omit the details since the proof is similar to the case of a Lie algebra,  see \cite{D}.

The $\fh_\ev$-module structure $\fh_\od$ is equivalent to the fact that $\lambda$ and $\mu$ are representations satisfying conditions (\ref{hod}). Here also we omit the details of the computations.

Let us consider the squaring on $\fh_\od$. We have (for all $\lambda \in \mathbb{K},$ for all $x_1\in \fg_\od$ and $a_1\in \fk_\od$):
\[
\begin{array}{lcl}
S_\fh(\lambda(x_1+a_1))&=&r(\lambda x_1,\lambda a_1)+s_\fg(\lambda x_1)+s_\fk(\lambda a_1)\\[3mm]
&=&\lambda^2 r( x_1,a_1)+\lambda^2 s_\fg(x_1)+\lambda^2 s_\fk(a_1)=\lambda^2 S_\fh(x_1+a_1),
\end{array}
\]

Besides, the map
\[
(x_1+a_1, y_1+b_1)\mapsto B(x_1+a_1,y_1+b_1):=S_\fh(x_1+y_1+a_1+b_1)-S_\fh(x_1+a_1)-S_\fh(y_1+b_1)
\]
is a~bilinear map on $\fh_\od$ with values in $\fh_\ev$. Indeed,
\[
\begin{array}{lcl}
 B(x_1+a_1, y_1+b_1)&=& r(x_1+y_1, a_1+b_1)+s_\fg(x_1+y_1) + s_\fk(a_1+b_1)+r(x_1, a_1) + s_\fg(x_1)\\[3mm]
&&+ s_\fk(a_1) + r(y_1, b_1) + s_\fg(y_1)-s_\fk(b_1)\\[3mm]
&=&r(x_1,b_1)+r(y_1,a_1)+[x_1,y_1 ]_\fg+[a_1,b_1 ]_\fk.
\end{array}
\]
Now let us consider the Jacobi identity involving to the squaring; for all $h=x_1+a_1 \in \fh_\od$ and $k\in \fh$ we have
\begin{itemize} \item The case where $k=y_0 + b_0\in \fh_\ev$. The LHS reads
\[
\begin{array}{lcl}
[S_\fh(h),k]_\fh&=&[r(x_1,a_1)+s_\fg(x_1)+s_\fk(a_1), y_0+b_0 ]_\fh + [r_\fg(x_1, a_1),y_0]_\fg\\[3mm]
&&+\rho (r_\fk(x_1, a_1))(y_0)+\pi(r_\fg(x_1, a_1))(b_0)+\pi(y)(r_\fk(x_1, a_1))+\rho(b)(r_\fg(x_1, a_1))\\[2mm]
&&+[r_k(x_1, a_1),b_0]_\fk +[x_1, [x_1,y_0]_\fg]_\fg+\pi(s_\fg(x_1))(b_0)+\rho(b)(s_\fg(x_1))+\rho (s_\fk(a_1))(y_0)\\[3mm]
&&+\pi(y)(s_\fk(a_1))+[a_1,[a_1,b_0]_\fk]_\fk.
\end{array}
\]
The RHS reads
\[
\begin{array}{lcl}
[h,[h,k]_\fh]_\fh&=&[x_1+a_1,[x_1,y_0]_\fg+\mu(b_0)(x_1)+\tilde \mu(b_0)(x_1)+\lambda(y_0)(a_1)+\tilde \lambda(y_0)(a_1)+[a_1, b_0]_\fk]_\fh\\[3mm]
&=&[x_1, [x_1, y_0]_\fg]_\fg+[x_1, \mu(b_0)(x_1)]_\fg+r(x_1, \tilde \mu(b)(x_1)) +r(x_1, \lambda(y)(a_1))\\[3mm]
&&+[x_1, \tilde \lambda(y)(a_1)]_\fg+r(a_1, [x_1, y]_\fg)+r(a_1, \mu(b)(x_1))+[a_1, \tilde \mu(b)(x_1)]_\fk\\[3mm]
&&+[a_1, \lambda(y)(a_1)]_\fk+r(a_1, \tilde \lambda(y)(a_1))+r(x_1, [a_1, b]_\fk)+[a_1, [a_1, b_0]_\fk]_\fk.
\end{array}
\]
Therfore, LHS=RHS is equivalent to conditions (\ref{sq5}) and (\ref{sq6}).
\item The case where $k=y_1 + b_1\in \fh_\od$. The LHS reads
\[
\begin{array}{lcl}
[S_\fh(h),k]_\fh&=&[r(x_1,a_1)+s_\fg(x_1)+s_\fk(a_1), y_1+b_1 ]_\fh\\[3mm]
&=&[r_\fg(x_1, a_1),y_1]_\fg++ \lambda(r_\fg(x_1, a_1), b)+\tilde \lambda (r_\fg(x_1, a_1), b)+\mu (r_\fk(x_1, a_1))(y_1) \\[3mm]
&&+\tilde \mu (r_\fk(x_1, a_1))(y_1)+[s_\fg(x_1),y_1]_\fg +\mu(s_\fk(a_1))(y_1)+\tilde \mu(s_\fk(a_1))(y_1)\\[3mm]
&&+[r_\fk(x_1, a_1), b_1]_\fk +\lambda (s_\fg(x_1))(b_1)+ \tilde \lambda (s_\fg(x_1))(b_1)+[s_\fk(a_1),b_1]_\fk.
\end{array}
\]
The RHS reads
\[
\begin{array}{lcl}
[h,[h,k]_\fh]_\fh&=&[x_1+a_1,[x_1, y_1]_\fg+r(y_1,a_1)+r(x_1, b_1)+[a_1, b_1]_\fk]_\fh\\[3mm]
&=&[x_1, [x_1, y_1]_\fg]_\fg+\lambda([x_1, y_1]_\fg)(a_1)+\tilde \lambda([x_1, y_1]_\fg)(a_1)+[a_1,[a_1,b]_\fk]_\fk \\[3mm]
&& + \mu([a,b]_\fk)(x_1)+\tilde \mu ([a,b]_\fk)(x_1)+\mu(r_\fk(x_1, a_1))(x_1)+\tilde \mu(r_\fk(x_1, a_1))(x_1)\\[3mm]
&&+[r_\fk(x_1, a_1),a_1]_\fk +[r_\fg(x_1,b),x_1]_\fg+\lambda(r_\fg(x_1,b))(a_1)+\tilde \lambda (r_\fg(x_1,b))(a_1)\\[3mm]
&&+\mu(r_\fk(a_1,y))(x_1)+\tilde \mu (r_\fk(a_1,y))(x_1)+[r_\fk(a_1, y),a_1]_\fk +[x_1, r_\fg(a_1,y)]_\fg \\[3mm]
&&+\lambda (r_\fg(a_1,y))(a_1)+\tilde \lambda (r_\fg(a_1,y))(a_1).
\end{array}
\]
Therfore, LHS=RHS is equivalent to conditions (\ref{sq7}) and (\ref{sq8}).\qed
\end{itemize} \noqed \end{proof}

\begin{proof}[Proof of Theorem \ref{mainmanin}] The proof will be a~consequence of Theorem \ref{mainthm}.  Let us explicitly describe the maps $\rho, \pi, \lambda, \tilde \lambda, \mu,$ and $\tilde \mu$:
\[\begin{array}{lcllcl}
\pi(x)(f)& = & f  \circ \ad_x \; \text{($x\in \fg_\ev$ and $f\in \fg_\ev^*$)}, &  \lambda(x)(f) & = &f \circ \ad_x \; \text{($x\in \fg_\ev$ and $f\in \fg_\od^*$),}\\[2mm]
\rho(f)(x) & = & x \circ \ad_f, \; \text{($x\in \fg_\ev$ and $f\in \fg_\ev^*$)}, &  \mu(f)(x)& = & x \circ \ad_f,\; \text{($x\in \fg_\od$ and $f\in \fg_\ev^*$),}\\[2mm]
\tilde \lambda(x)(f) & = & x \circ \ad_f, \text{($x\in \fg_\ev$ and $f\in \fg_\od^*$)}, &  \tilde \mu (f)(x)& = & f\circ \ad_x, \;\text{($x\in \fg_\ev$ and $f\in \fg_\od$).}
\end{array}
\]
On the other hand, the maps $r_\fg$ and $r_{\fg^*}$ are given by (where $f\in \fg_\od^*$ and $x\in \fg_\od$)
\[
r_\fk(x,f)=f\circ \ad_x, \quad r_\fg(x,f)=x\circ \ad_f.
\]
It follows that the squaring on $\fh$ is given by
\[
s_\fh(x+f)=s_\fg(x)+s_{\fg^*}(f)+f\circ \ad_x+ x\circ \ad_f, \; \text{for all $x\in \fg_\od$ and $f\in \fg_\od^*$}.
\]
A direct computation shows that the conditions (\ref{Sq}), (\ref{Sq*}) and (\ref{Bra}) are equivalent to (\ref{hev}), (\ref{hod}), (\ref{sq5}), (\ref{sq6}), (\ref{sq7}) and  (\ref{sq8}). 
The proof of Theorem \ref{mainmanin} is now complete.
\end{proof}

\section{$\od$-antisymmetric forms,  Manin Triples and double extensions}\label{symplectic}

\subsection{Double extensions (DEs) of even NIS-Lie superalgebras equipped with closed $\od$-antisymmetric ortho-orthogonal forms}\label{DEeven}
The purpose of this subsection is to study a class of Lie superalgebras which admit an admissible  Yang-Baxter operator $U$ as in Prop.  \ref{propU}.  As a result, such a Lie superalgebra from this class carries a Lie bi-superalgebra structure, and hence it is  a Manin triple.  Throughout the process of {\it double extensions}, the operator $U$ is constructed inductively.  Double extensions of Lie superalgebras in characteristic 2 were introduced and studied in \cite{BeBou}. For interesting examples of DEs of simple Lie superalgebras classified up to isomorphisms, see \cite{BeBou, BLS}. 

Before introducing the process of double extensions, we will give a proposition that makes a link between homogenous invertible $\mathscr B$-symmetric derivations and admissible Yang-Baxter operators.  We will need the following lemma. \\

\begin{mLemma}\label{Uinv}
Let $\fg$ be Lie superalgebra in characteristic 2, and let $U \in \Aut(\fg)$ be an even automorphism. Then $U^{-1}$ is a derivation on $\fg$ if and only if 

\textup{(i)} $U([U(x),y]_\fg+[x, U(y)]_\fg)=[U(x), U(y)]_\fg$ for all $x,y \in \fg$.

\textup{(ii)} $U([x,U(x)])=s_\fg(U(x))$ for all $x\in \fg_\od$.
\end{mLemma}
\begin{proof} Straightforward. 
\end{proof}
\begin{mProposition}\label{Inverse0}
Let $\fg$ be Lie superalgebra  in characteristic 2 with an~even NIS $\mathscr B$.  Let $\Delta$ be an even $\mathscr B$-symmetric invertible derivation. Then $U=\Delta^{-1}$ is an admissible Yang-Baxter operator on $\fg.$
\end{mProposition}
\begin{proof}
Apply Lemma \ref{Uinv} then Proposition \ref{propU}.
\end{proof}
{\bf Example 4.} Consider the Lie superalgebra $\fg=\fh\oplus \fh^*$ as in Example 3 of \S \ref{examples}.   The map $U$ defined by 
\[
U(x+f):={\mathscr D}^{-1}(x)+f\circ {\mathscr D}^{-1} \quad \text{for all $x\in \fh$ and $\fh^*$},
\]
is an admissible Yang-Baxter operator on $\fg$.\\

The first two theorems below have been proved in \cite{BeBou}.\\

\begin{mTheorem}[${\mathscr D}_\ev$-extension -- NIS even]\label{MainTh} Let $(\fa, {\mathscr B}_\fa)$ be a~NIS-Lie superalgebra in~charac-teristic $2$ such that ${\mathscr B}_\fa$ is even. Let ${\mathscr D}\in \fder_\ev(\fa)$ be a~derivation satisfying the following conditions:
\begin{eqnarray}
\label{D1} {\mathscr B}_\fa({\mathscr D}(a),b)+{\mathscr B}_\fa(a,{\mathscr D}(b))=0\text{ for any } a,b \in \fa; \; {\mathscr B}_\fa({\mathscr D}(a), a)=0\text{ for any } a~\in \fa_\ev.
\end{eqnarray}
Let $\alpha: \fa_\od \rightarrow \Bbb K$ be a~quadratic form and ${\mathscr B}_\alpha$ its polar form\footnote{In characteristic $2$, the polar form ${\mathscr B}_\alpha$ associated with a~quadratic form $\alpha$ is defined by \[(a,b)\mapsto {\mathscr B}_\alpha(a,b):=\alpha(a+b)+\alpha(a)+\alpha(b).\]} which satisfies 
\begin{eqnarray}
\label{D3}
{\mathscr B}_\fa(a, {\mathscr D}(b))&=&{\mathscr B}_\alpha(a,b) \text{ for any  $a,b \in \fa_\od$}. 
\end{eqnarray}

Then,  there exists a~NIS-Lie superalgebra structure on $\fg:=\mathscr{K} \oplus \fa \oplus \mathscr{K} ^*$, where $\mathscr{K} :=\Span\{x\}$ for $x$ even, defined as follows. The squaring is given by 
\[
s_\fg(a):= s_\fa(a)+\alpha (a) x \qquad \text{ for any } a\in \fg_\od \; (=\fa_\od).
\]
The bracket on $\fg$ is defined as follows:
\be\label{*}
[\mathscr{K} ,\fg]_\fg=0, \qquad [a,b]_\fg:={\mathscr B}_\fa({\mathscr D}(a),b)x+[a,b]_\fa, \qquad [x^*,a]_\fg:={\mathscr D}(a)\text{ for any } a,b\in  \fa.
\ee
The non-degenerate $\ev$-antisymmetric bilinear form ${\mathscr B}_\fg$ on $\fg$ is defined as follows: 
\begin{eqnarray*}
&{{\mathscr B}_\fg}\vert_{\fa \times \fa}:= {\mathscr B}_\fa, \quad {\mathscr B}_\fg(\fa,\mathscr{K} ):=0, \quad {\mathscr B}_\fg(x,x^*):=1, \quad {\mathscr B}_\fg(\fa,\mathscr{K} ^*):=0, \\ 
&{\mathscr B}_\fg(x,x):=0,\ {\mathscr B}_\fg(x^*,x^*) \; \text{ arbitrary}.
\end{eqnarray*}
Moreover, the form ${\mathscr B}_\fg$, obviously even, is invariant on $\fg$. 

Therefore, $(\fg,{\mathscr B}_\fg)$ is a~NIS-Lie superalgebra. \\

\end{mTheorem}

\begin{mTheorem}[${\mathscr D}_\od$-extension -- NIS even]\label{MainTh2} Let $(\fa,{\mathscr B}_\fa)$ be a~NIS-Lie superalgebra in charac-teristic $2$ such that ${\mathscr B}_\fa$ is even. Let ${\mathscr D}\in\fder_\od(\fa)$ and $a_0\in \fa_\ev$ satisfy the following conditions:
\begin{eqnarray}
\label{2D1} {\mathscr B}_\fa(\mathscr{D}(a),b)&=& {\mathscr B}_\fa(a,\mathscr{D}(b)) \text{ for any } a,b \in\fa;\\[2mm]
\label {2D2} {\mathscr D}^2&=&\ad_{a_0};\\[2mm]
\label{2D3} {\mathscr D}(a_0)&=&0.
\end{eqnarray}
Then there exists a~NIS-Lie superalgebra structure on $\fg:=\mathscr{K} \oplus \fa \oplus \mathscr{K}^*$, where $\mathscr{K}:=~\Span\{x\}$ and $x$ odd, defined as follows. The squaring is given by 
\[
s_\fg(rx+a+t x^*):= s_\fa(a) + t^2 a_0+t {\mathscr D}(a) \qquad \text{ for any }\; rx+a+tx^* \in \fg_\od.
\]
The bracket is given by:
\[
\begin{array}{l}
[x,\fg]_\fg:=0; \qquad [a,b]_\fg:=[a,b]_\fa+{\mathscr B}_\fa(\mathscr{D}(a),b)x\quad\text{~~for any $a,b\in\fa$};\\[2mm]
[x^*,a]_\fg:={\mathscr D}(a) \quad  \text{ for any } a~\in \fa.
\end{array}
\]
The bilinear form ${\mathscr B}_\fg$ on $\fg$ defined by:
\begin{eqnarray*}
&{{\mathscr B}_\fg}\vert_{\fa \times \fa}:= {\mathscr B}_\fa, \quad {\mathscr B}_\fg(\fa,\mathscr{K}):=0, \quad {\mathscr B}_\fg(\fa,\mathscr{K}^*):=0, \\ 
&{\mathscr B}_\fg(x,x^*):=1, \quad {\mathscr B}_\fg(x,x):= {\mathscr B}_\fg(x^*,x^*):=0.
\end{eqnarray*}
is even,  non-degenerate,  $\ev$-antisymmetric and invariant on $(\fg, [\cdot ,\cdot]_\fg, s_\fg)$.

Therefore, $(\fg,{\mathscr B}_\fg)$ is a~NIS-Lie superalgebra. \\

\end{mTheorem}
%
In what follows, we will consider NIS Lie superalgebras endowed with a closed $\od$-antisymmetric ortho-orthogonal form.  Example 3 in \S \ref{examples} provides an example of such Lie superalgebras. \\

\begin{mTheorem}[Ortho-orthogonal forms for ${\mathscr D}_\ev$-extensions] \label{MainDer} Let $(\fa, {\mathscr B}_\fa, \omega_\fa)$ be an~even NIS Lie superalgebra equipped with a closed $\od$-antisymmetric ortho-orthogonal form $\omega_\fa$,  and let $\tilde \Delta\in \fder(\fa)_\ev$ be the unique and invertible derivation such that $\omega_\fa(\cdot,\cdot)={\mathscr B}_\fa(\tilde \Delta (\cdot),\cdot)$. Let $(\fg, {\mathscr B}_\fg)$ be the double extension of $(\fa, {\mathscr B}_\fa)$ by means of a~derivation ${\mathscr D}\in \fder(\fg)_\ev$ as in Theorem \ref{MainTh} such that ${\mathscr B}_\fg(x^*,x^*)=0$. Suppose there exist $\lambda\in \Kee \backslash \{0\}$ and $a_0\in \fa_\ev$ such that 
\begin{eqnarray}
\label{eq7} [\tilde \Delta, {\mathscr D}]& = &\lambda {\mathscr D}+\ad_{a_0},\\
\label{eq8}\lambda \alpha(a)+{\mathscr B}_{\alpha}({\mathscr D}\circ \tilde \Delta (a),a)& = & {\mathscr B}_\fa(s_\fa(a),a_0) \text{ for all $a\in \fa_\od$}.
\end{eqnarray}
The linear map $\Delta: \fg \rightarrow \fg$ defined by
\[
\begin{array}{lcl}
\Delta(x)&=& \lambda x,\\
\Delta(x^*)&=&\lambda x^*+a_0,\\
\Delta(a)&=& \tilde \Delta(a)+{\mathscr B}_\fa(a,a_0)x.
\end{array}
\]
is an invertible derivation on $\fg$ for which ${\mathscr B}_\fg$ is $\Delta$-invariant. Consequently, $(\fg, {\mathscr B}_\fg, \omega_\fg)$ is a~NIS-Lie superalgebra with a closed $\od$-antisymmetric ortho-orthogonal form $\omega_\fg:={\mathscr B}_\fg(\Delta(\cdot),\cdot)$.

\end{mTheorem}
By calculating directly and applying Prop. \ref{Inverse0},  we get the following Corollary. \\

\begin{mCorollary} Let us consider the admissible Yang-Baxter operator $\tilde U:=\tilde \Delta^{-1}$ on $\fa$. The map $U:=\Delta^{-1}$ where
\[
\begin{array}{lcl}
U(x)&=& \lambda^{-1} x,\\
U(x^*)&=&\lambda^{-1} x^*+\lambda^{-1} \tilde U(a_0)+ \lambda^{-2}{\mathscr B}_\fa(\tilde U(a_0), a_0)x,\\[1mm]
U(a)&=& \tilde U(a)+ \lambda^{-1} {\mathscr B}_\fa(\tilde U(a),a_0)x, 
\end{array}
\]
is an admissible Yang-Baxter operator on the double extension $\fg$.
\end{mCorollary}
\begin{proof} [Proof of Theorem \ref{MainDer}]
First, let us observe that $\Delta$ is obviously invertible since $\tilde \Delta$ is invertible. Let us prove that $\Delta$ is a~derivation. We will check the condition for the squaring. Indeed,  for any $a\in \fa_\od$
\[
\Delta(s_\fg(a))= \Delta(s_\fa(a)+\alpha(a)x)=\Delta(s_\fa(a))+ \lambda \alpha(a)x=\tilde \Delta(s_\fa(a))+{\mathscr B}_\fa(s_\fa(a),a_0)x+ \lambda \alpha(a)x.
\]
Aditionally, we have for any $a\in \fa_\od$
\[
[\Delta(a),a]_\fg=[\tilde \Delta(a),a]_\fg=[\tilde \Delta(a),a]_\fa+{\mathscr B}_\alpha({\mathscr D}\circ \tilde \Delta(a),a)x.
\]
The fact that $\tilde \Delta \in \fder(\fa)_\ev$ and Eq. (\ref{eq8}) imply that $\Delta(s_\fg(a))=[\Delta(a),a]_\fg$. 

Now let us check the condition
\[
\Delta([f,g]_\fg)=[\Delta(f),g]_\fg+[f, \Delta(g)]_\fg \text{ for all $f,g \in \fg$}.
\]
We will only check the case where $(f,g)=(a,x^*)$; checking the other cases is a~routine. Indeed, 
\[
\begin{array}{lcl}
\Delta([a,x^*]_\fg)&=&\Delta({\mathscr D}(a))= \tilde \Delta({\mathscr D}(a))+{\mathscr B}_\fa({\mathscr D}(a),a_0)x.
\end{array}
\]
On the other hand, 
\[
\begin{array}{lcl}
[\Delta(a),x^*]_\fg+[a, \Delta(x^*)]_\fg & = &[\tilde \Delta(a),x^*]_\fg+[a,\lambda x^*+a_0]_\fg\\[2mm]
&=& {\mathscr D}(\tilde \Delta(a))+\lambda {\mathscr D}(a) + [a,a_0]_\fa+ {\mathscr B}_\fa({\mathscr D}(a), a_0)x.
\end{array}
\]
Therefore,  $\Delta([a,x^*]_\fg)=[\Delta(a),x^*]_\fg+[a, \Delta(x^*)]_\fg$ is equivalent to Eq. (\ref{eq7}). 

Let us show that ${\mathscr B}_\fg$ is $\Delta$-invariant. Indeed, for all $a\in \fa_\ev$ (recall that $x$ is orthogonal to $\fa$):
\[
{\mathscr B}_\fg(\Delta(a),a)={\mathscr B}_\fg(\tilde \Delta(a)+{\mathscr B}_\fa(a,a_0)x,a)={\mathscr B}_\fa(\tilde \Delta(a),a)=\omega_\fa(a,a)=0. 
\]
Similarly, 
\[
{\mathscr B}_\fg(\Delta(x^*),x^*)={\mathscr B}_\fg(\lambda x^*+a_0, x^*)=0 \text{ and } {\mathscr B}_\fg(\Delta(x),x)={\mathscr B}_\fg(\lambda x,x)=0.
\]
For all $a,b\in \fa$, we have
\[
{\mathscr B}_\fg(\Delta(a),b)={\mathscr B}_\fg(\tilde \Delta(a) + {\mathscr B}_\fa(a,a_0)x,b)={\mathscr B}_\fa(\tilde \Delta(a),b)={\mathscr B}_\fa(a,\tilde \Delta(b))={\mathscr B}_\fg(a,\Delta(b)).
\]
Moreover, 
\[
\begin{array}{lcl}
{\mathscr B}_\fg(\Delta(x^*),a) + {\mathscr B}_\fg(x^*, \Delta(a))& = & {\mathscr B}_\fg(\lambda x^*+a_0, a)+ {\mathscr B}_\fg(x^*, \tilde \Delta(a)+{\mathscr B}_\fa(a,a_0)x.) \\[2mm]
&= &{\mathscr B}_\fa(a_0,a)+{\mathscr B}_\fa(a,a_0)=0.
\end{array}
\]
The $\Delta$-invariant of  ${\mathscr B}_\fg$ can easily be checked for the remaining elements of the basis of $\fg$.\end{proof}

In proving the converse of Theorem \ref{MainTh} in \cite{BeBou}, the following object, called special center, has been introduced
\[
\fz_s(\fg):=\fz(\fg) \cap s_\fg(\fg_\od)^\perp.
\]
\begin{mLemma}[Description of the special centre]
If the Lie superalgebra $(\fg, {\mathscr B})$ admits an invertible $\Delta$ derivation such that $ {\mathscr B}$ is $\Delta$-invariant, then $\fz_s(\fg)=\fz(\fg)$.

\end{mLemma}
\begin{proof}

Let $x\in \fz(\fg)$. There exists $\tilde x \in \fz(\fg)$ such that $x=\Delta(\tilde x)$. We have
\[
{\mathscr B}_\fg(x, s_\fg(y))= {\mathscr B}_\fg(\Delta(\tilde x), s_\fg(y))={\mathscr B}_\fg(\tilde x, \Delta(s_\fg(y)))={\mathscr B}_\fg(\tilde x, [y,\Delta(y)]_\fg)={\mathscr B}_\fg([\tilde x, y]_\fg,\Delta(y))=0.
\]
\end{proof}
In Theorem \ref{Rec1} below, the ground field $\mathbb K$ is assumed to be algebraically closed. \\
\begin{mTheorem}[Converse of Theorem \ref{MainDer}]
\label{Rec1}
Let $(\fg,\mathscr{B}_\fg, \omega_\fg)$ be an irreducible even NIS-Lie superalgebra equipped with a closed $\od$-antisymmetric ortho-orthogonal form $\omega_\fg$. Suppose that $\fz_s(\fg)_\ev\not=\{0\}.$ Then,  $(\fg,\mathscr{B}_\fg, \omega_\fg)$ is obtained as a $\mathscr{D}$-extension of a~NIS-Lie superalgebra $(\fa,\mathscr{B}_\fa, \omega_\fa)$ with a closed $\od$-antisymmetric ortho-orthogonal form $\omega_\fa$.

\end{mTheorem}
\begin{proof} 

Let $\Delta\in \fder(\fg)_\ev$ be the unique invertible derivation such that $\omega_\fg(\cdot,\cdot)=\mathscr{B}_\fg(\Delta(\cdot),\cdot)$. It follows that $\Delta(\fz_s(\fg)_\ev)=\fz_s(\fg)_\ev$. Since $\Kee$ is algebraically closed, let us choose $x\in \fz_s(\fg)_\ev$ such that $\Delta(x)=\lambda x$ for some $\lambda \in \Kee \backslash \{0\}$. It has been shown in \cite[Prop. 3.1.3]{BeBou}, that the subspaces $\mathscr{K}:=\Span\{x\}$ and $\mathscr{K}^\perp$ are ideals in $(\fg,{\mathscr B}_\fg)$. Besides, $\Delta(\mathscr{K}^\perp) \subset \mathscr{K}^\perp$. Indeed,  for all $f\in \mathscr{K}^\perp$, we have 
\[
{\mathscr B}_\fg(\Delta(f),x)={\mathscr B}_\fg(f,\Delta(x))={\mathscr B}_\fg(f,\lambda x)=0.
\]
Now, since $\mathscr{K}$ is 1-dimensional, then either $\mathscr{K} \cap \mathscr{K}^\perp=\{0\}$ or $\mathscr{K}\cap \mathscr{K}^\perp=\mathscr{K}$. The first case is to be disregarded because otherwise $\fg=\mathscr{K}\oplus \mathscr{K}^\perp$ and the Lie superalgebra $\fg$ is not irreducible. Hence, $\mathscr{K}\cap \mathscr{K}^\perp=\mathscr{K}$. It follows that $\mathscr{K}\subset \mathscr{K}^\perp$ and $\dim(\mathscr{K}^\perp)=\dim(\fg)-1$. Therefore, there exists a non-zero $x^* \in \fg_\ev$ such that
\[
\fg=\mathscr{K}^\perp\oplus \mathscr{K}^*, \quad \text{ where  $\mathscr{K}^*:=\Span\{x^*\}$}.
\]
This $x^*$ can be normalized such that ${\mathscr B}_\fg(x,x^*)=1$. Besides, ${\mathscr B}_\fg(x,x)=0$ since $\mathscr{K}\cap \mathscr{K}^\perp=\mathscr{K}$.

Let us define $\fa:=(\mathscr{K} +\mathscr{K}^*)^\perp$. We then have a~decomposition $\fg=\mathscr{K} \oplus \fa \oplus \mathscr{K}^*$. 

It has been shown in \cite[Prop. 3.1.3]{BeBou} that $(\fg, {\mathscr B}_\fg)$ is a~double extension of $(\fa, {\mathscr B}_\fa)$, where $ {\mathscr B}_\fa= {\mathscr B}_{\fg_{|_{\fa\times \fa}}}$, by means of a~derivation $\mathscr D$ and a~quadratic form $\alpha$.

Let us write $\Delta$ as ($a_0\in \fa_\ev, \beta(a), \gamma(a), \mu, \nu \in \Kee$):
\[
\begin{array}{lcl}
\Delta(x)&=& \lambda x,\\
\Delta(x^*)&=&\nu x^*+a_0+\mu x,\\
\Delta(a)&=& \tilde \Delta(a)+\beta(a)x + \gamma(a) x^*.
\end{array}
\]
The condition ${\mathscr B}_\fg(\Delta(x^*),x)={\mathscr B}_\fg(x^*, \Delta(x))$ implies that $\nu=\lambda$. Similarly, the condition ${\mathscr B}_\fg(\Delta(x^*),x^*)=0$ implies that $\mu=0$ and the condition ${\mathscr B}_\fa(\Delta(a),x)={\mathscr B}_\fa(a,\Delta(x))$ implies that $\gamma(a)=0$ for every $a\in \fa$.

In addition,  the conditions 
\[
\begin{array}{lcl}
{\mathscr B}_\fg(\Delta(a),x^*) & = & {\mathscr B}_\fg ( \tilde \Delta(a)+\beta(a)x,x^*)=\beta(a), \\[1mm]
{\mathscr B}_\fg(a, \Delta(x^*)) & = & {\mathscr B}_\fg(a, \lambda x^*+a_0+\mu x)={\mathscr B}_\fa(a,a_0),
\end{array}
\]
imply that 
\begin{equation}\label{eq1}
\beta(a)={\mathscr B}_\fa(a,a_0) \text{ for all $a\in \fa$}.
\end{equation}
Now, we have 
\[
\begin{array}{lcl}
\Delta([x^*,a]_\fg)&=&\Delta({\mathscr D}(a))=\tilde \Delta ({\mathscr D}(a))+\beta({\mathscr D}(a))x\; \text{ for all $a\in \fa$}.
\end{array}
\]
On the other hand, we have
\[
[\Delta(x^*),a]_\fg+[x^*,\Delta(a)]_\fg=\lambda {\mathscr D}(a)+[a_0,a]_\fa+{\mathscr B}_\fa({\mathscr D}(a_0),a)x+{\mathscr D}(\tilde \Delta(a)) \text{ for all $a\in \fa$}.
\]
Since $\Delta$ is a~derivation on $\fg$, it follows that 
\begin{equation}
\label{eq2}
[\tilde \Delta, {\mathscr D}]=\lambda {\mathscr D}+\ad_{a_0}. 
\end{equation}
Following the same computation as in \cite{BaB}, the condition 
\[
\Delta([a,b]_\fg)=[\Delta(a),b]_\fg+[a,\Delta(b)]_\fg \; \text{ for all $a,b\in \fa$,}
\]
implies that 
\begin{eqnarray}
\label{eq3} \tilde \Delta([a,b]_\fa)&=&[\tilde \Delta(a),b]_\fa+[a, \tilde \Delta(b)]_\fa \\
\label{eq4} \beta([a,b]_\fa)&=&{\mathscr B}_\fa({\mathscr D}(a),b)+ {\mathscr B}_\fa({\mathscr D}(\tilde \Delta(a)),b)+ {\mathscr B}_\fa({\mathscr D}(a), \tilde \Delta(b)).
\end{eqnarray}
Eq. (\ref{eq4}) is certainly satisfied, thanks to Eqs. (\ref{eq1}), (\ref{eq2}) and the fact that ${\mathscr B}_\fa$ is $\tilde \Delta$-invariant.

For all $a\in \fa$, the condition $\Delta(s_\fg(a))=[a,\Delta(a)]_\fg$ implies that 
\begin{eqnarray}
\label{eq5pre} \tilde \Delta (s_\fa(a)) & = & [a,\tilde \Delta (a)]_\fa,\\
\label{eq6} {\mathscr B}_\fa(s_\fa(a),a_0)&=&\lambda \alpha(a)+{\mathscr B}_{\alpha}(a, \tilde \Delta (a)).
\end{eqnarray}
Eqs. (\ref{eq3}), (\ref{eq5pre}) imply that $\tilde \Delta$ is a~derivation on $\fa$. Since $\Delta$ is invertible, then $\tilde \Delta$ is invertible as well.  Define
\[
\omega_\fa(a,b):={\mathscr B}_\fa(\tilde \Delta (a), b) \; \text{ for all $a,b\in \fa$}.
\]
The form $\omega_\fa$ is closed,  $\od$-antisymmetric and ortho-orthogonal over $\fa$.\end{proof}

Hereafter, we will carry out the same study for the $\mathscr D_\od$-extensions.\\

\begin{mTheorem}[Ortho-orthogonal forms for ${\mathscr D}_\od$-extensions] \label{MainDer2} Let $(\fa, {\mathscr B}_\fa, \omega_\fa)$ be an~even NIS Lie superalgebra equipped with a closed $\od$-antisymmetric ortho-orthogonal form $\omega_\fa$, and let $\tilde \Delta\in \fder(\fg)_\ev$ be the unique and invertible derivation such that $\omega_\fa(\cdot,\cdot)={\mathscr B}_\fa(\tilde \Delta (\cdot),\cdot)$. Let $(\fg, {\mathscr B}_\fg)$ be the double extension of $(\fa, {\mathscr B}_\fa)$ by means of a~derivation ${\mathscr D}\in \fder_\od(\fa)$ satisfying \eqref{2D1}, \eqref{2D2}, \eqref{2D3}. Suppose there exist $\lambda\in \Kee \backslash \{0\}$ and $b_0\in \fa_\od$ such that 
\begin{eqnarray}
\label{eqb9} [\tilde \Delta, {\mathscr D}]& = &\lambda {\mathscr D}+\ad_{b_0},\\
\label{eqb10} \tilde \Delta(a_0)& = & {\mathscr D}(b_0).
\end{eqnarray}
The linear map $\Delta: \fg \rightarrow \fg$ defined by \textup{(}$\mu\in \Kee$ is arbitrary\textup{)}
\[
\begin{array}{lcl}
\Delta(x)&=& \lambda x,\\
\Delta(x^*)&=&\lambda x^*+b_0+\mu x,\\
\Delta(a)&=& \tilde \Delta(a)+{\mathscr B}_\fa(a,b_0)x.
\end{array}
\]
is an invertible derivation on $\fg$ such that ${\mathscr B}_\fg$ is $\Delta$-invariant. Consequently, $(\fg, {\mathscr B}_\fg, \omega_\fg)$ is a~NIS-Lie superalgebra with a closed $\od$-antisymmetric ortho-orthogonal form $\omega_\fg:={\mathscr B}_\fg(\Delta(\cdot),\cdot)$.\\

\end{mTheorem}

By applying Prop. \ref{Inverse0},  we get the following Corollary. \\

\begin{mCorollary}
Let us consider the admissible Yang-Baxter operator $\tilde U:=\tilde \Delta^{-1}$ on $\fa$. The map $U:=\Delta^{-1}$ where
\[
\begin{array}{lcl}
U(x)&=& \lambda^{-1} x,\\
U(x^*)&=&\lambda^{-1} x^*+ \lambda^{-1} \tilde U(b_0)+\lambda^{-1}(\lambda^{-1} {\mathscr B}_\fa (\tilde U(b_0), b_0)+\mu x),\\
U(a)&=& \tilde U(a)+\lambda^{-1} {\mathscr B}_\fa(\tilde U(a),b_0)x,
\end{array}
\]
is an admissible Yang-Baxter operator on the double extension $\fg$.
\end{mCorollary}
\begin{proof} [Proof of Theorem \ref{MainDer2}]

Obviously,  the map $\Delta$ is invertible since $\tilde \Delta$ is invertible. Let us check that $\Delta$ is a~derivation $\fg$. We start by checking the condition on the squaring. Indeed,  for all $a\in \fa_\od$ and $s,t\in \Kee$
\[
\begin{array}{l}
\Delta(s_\fg(rx+a+tx^*)) \\[2mm]
=\Delta( s_\fa(a) + t^2 a_0+t {\mathscr D}(a))\\[2mm]
= \tilde \Delta(s_\fa(a) + t^2 a_0+t {\mathscr D}(a))+{\mathscr B}_\fa(s_\fa(a) + t^2 a_0+t {\mathscr D}(a),b_0)x\\[2mm]
= \tilde \Delta(s_\fa(a)) + t^2 \tilde \Delta(a_0)+t \tilde \Delta({\mathscr D}(a)).
\end{array}
\]
Besides, 
\[
\begin{array}{l}
[\Delta(rx+a+tx^*),rx+a+tx^*]_\fg\\[2mm]
 =  [r \lambda x+\tilde \Delta(a)+{\mathscr B}_\fa(a,b_0)x+t(\lambda x^*+b_0+\mu x), rx+a+tx^*]_\fg\\[2mm]
 =  [\tilde \Delta(a)+t(\lambda x^*+b_0), a+tx^*]_\fg\\[2mm]
 =  [\tilde \Delta(a),a]_\fa+t  {\mathscr D}(\tilde \Delta(a))+t \lambda {\mathscr D}(a)+t [b_0,a]_\fa+t^2 {\mathscr D}(b_0).
\end{array}
\]
The fact that $\tilde \Delta \in \fder(\fa)_\od$ and Eqs. (\ref{eqb9}), (\ref{eqb10}) imply that $\Delta(s_\fg(f))=[\Delta(f),f]_\fg$ for all $f\in \fg_\od$. 

Now let us check the condition
\[
\Delta([f,g]_\fg)=[\Delta(f),g]_\fg+[f, \Delta(g)]_\fg \text{ for all $f,g \in \fg$}.
\]
We will only check the case where $(f,g)=(a,x^*)$; checking the other cases is routine. Indeed,  using Eq. (\ref{eqb9}) we have
\[
\begin{array}{lcl}
\Delta([a,x^*]_\fg)&=&\Delta({\mathscr D}(a))= \tilde \Delta({\mathscr D}(a))+{\mathscr B}_\fa({\mathscr D}(a),b_0)x.
\end{array}
\]
On the other hand, 
\[
\begin{array}{lcl}
[\Delta(a),x^*]_\fg+[a, \Delta(x^*)]_\fg & = &[\tilde \Delta(a)+{\mathscr B}(a,b_0)x,x^*]_\fg+[a,\lambda x^*+b_0+\mu x]_\fg\\[2mm]
&=& {\mathscr D}(\tilde \Delta(a))+\lambda {\mathscr D}(a) + [a,b_0]_\fa+ {\mathscr B}_\fa({\mathscr D}(a), b_0)x\\[2mm]
\end{array}
\]
Therefore, the condition 
\[
\Delta([a,x^*]_\fg)=[\Delta(a),x^*]_\fg+[a, \Delta(x^*)]_\fg
\] 
is equivalent to Eq. (\ref{eqb9}). 

Let us show that ${\mathscr B}_\fg$ is $\Delta$-invariant. Indeed, for all $a\in \fa_\ev$ (recall that $x$ is orthogonal to $\fa$):
\[
{\mathscr B}_\fg(\Delta(a),a)={\mathscr B}_\fg(\tilde \Delta(a)+{\mathscr B}_\fa(a,b_0)x,a)={\mathscr B}_\fa(\tilde \Delta(a),a)=0. 
\]
For all $a,b\in \fa$, we have
\[
{\mathscr B}_\fg(\Delta(a),b)={\mathscr B}_\fg(\tilde \Delta(a) + {\mathscr B}_\fa(a,b_0)x,b)={\mathscr B}_\fa(\tilde \Delta(a),b)={\mathscr B}_\fa(a,\tilde \Delta(b))={\mathscr B}_\fg(a,\Delta(b)).
\]
Similarly, 
\[
\begin{array}{lcl}
{\mathscr B}_\fg(\Delta(x^*),a) + {\mathscr B}_\fg(x^*, \Delta(a))& = & {\mathscr B}_\fg(\lambda x^*+b_0+\mu x, a)+ {\mathscr B}_\fg(x^*, \tilde \Delta(a)+{\mathscr B}_\fa(a,b_0)x) \\[2mm]
&= &{\mathscr B}_\fa(b_0,a)+{\mathscr B}_\fa(a,b_0)=0.
\end{array}
\]
The $\Delta$-invariance of ${\mathscr B}_\fg$ can easily be checked for the remaining elements of the basis of $\fg$. We define a form on $\fg$ as follows
\[
\omega(\cdot, \cdot):={\mathscr B}_\fg(\Delta(\cdot), \cdot).\qed
\]
This form is closed, $\od$-antisymmetric and ortho-orthogonal on $\fg$.
\noqed \end{proof}

In proving the converse of Theorem \ref{MainTh2} in \cite{BeBou},  we had to introduce the cone
\[
\mathscr{C}(\fg, {\mathscr B}_\fg):=\{x\in \fg_\od\; | \; {\mathscr B}_\fg(s_\fg(x), s_\fg(t))=0 \quad \text{ for any } t\in \fg_\od \},
\]
Clearly, if $x\in\mathscr{C}(\fg, {\mathscr B}_\fg)$, then $\lambda x\in\mathscr{C}(\fg, {\mathscr B}_\fg)$ for any $\lambda\in\Kee$, hence the terminology. \\

\begin{mLemma}[Description of the cone]
If the Lie superalgebra $(\fg, {\mathscr B}_\fg)$ admits an invertible derivation $\Delta$ such that  ${\mathscr B}_\fg$ is $\Delta$-invariant,  then 
\[
\fz(\fg)_\od \cap~\mathscr{C} (\fg, {\mathscr B}_\fg)=\fz(\fg)_\od.
\]

\end{mLemma}
\begin{proof}

Let $x\in \fz(\fg)_\od$. For every $t\in \fg_\od$, there exists a $\tilde t \in \fg_\ev$ such that $\Delta(\tilde t)=s_\fg(t)$. We have
\[
{\mathscr B}_\fg(s_\fg(x), s_\fg(t))={\mathscr B}_\fg(s_\fg(x), \Delta(\tilde t))={\mathscr B}_\fg(\Delta(s_\fg(x)), \tilde t)={\mathscr B}_\fg([x,\Delta(x)]_\fg, \tilde t)=0.
\]
Therefore, $x\in \mathscr{C} (\fg, {\mathscr B}_\fg)$ and we are done.
\end{proof}
In Theorem \ref{Rec2} below, the ground field $\mathbb K$ is assumed to be algebraically closed. \\
\begin{mTheorem}[Converse of Theorem \ref{MainDer2}]
\label{Rec2} Let $(\fg,\mathscr{B}_\fg, \omega_\fg)$ be an irreducible even NIS-Lie superalgebra and equipped with a closed $\od$-antisymmetric ortho-orthogonal form $\omega_\fg$. Suppose that  $\fz(\fg)_\od \not=\{0\}$. Then $(\fg,\mathscr{B}_\fg, \omega_\fg)$ is obtained as an $\mathscr{D}$-extension from a~NIS-Lie superalgebra $(\fa,\mathscr{B}_\fa, \omega_\fa)$,  where $\omega_\fa$ is closed, $\od$-antisymmetric and ortho-orthogonal form on~$\fa$.

\end{mTheorem}
\begin{proof}  

Let $\Delta\in \fder(\fg)_\ev$ be the unique invertible derivation such that $\omega_\fg(\cdot,\cdot)=\mathscr{B}_\fg(\Delta(\cdot),\cdot)$. It follows that $\Delta(\fz(\fg)_\od)=\fz(\fg)_\od$. Since $\Kee$ is algebraically closed, let us choose $x\in \fz(\fg)_\od$ such that $\Delta(x)=\lambda x$ for some $\lambda \in \Kee \backslash \{0\}$. We assume that $s_\fg(x)=0$; otherwise, $s_\fg(x)$ will be a~non-zero element of $\fz_\ev(\fg)$ and hence Theorem \ref{Rec1} can be applied.  Besides, $\Delta(\mathscr{K}^\perp) \subset \mathscr{K}^\perp$. Indeed, (for all $f\in \mathscr{K}^\perp$)
\[
{\mathscr B}_\fg(\Delta(f),x)={\mathscr B}_\fg(f,\Delta(x))={\mathscr B}_\fg(f,\lambda x)=0.
\]
It has been shown in \cite[Prop. 3.2.2]{BeBou} that the subspace $\mathscr{K}:= \Span\{x\}$ is an ideal in $\fg$; moreover, the subspace $\mathscr{K}^\perp$ is also an ideal containing $\mathscr{K}$. Same arguments as in Theorem \ref{Rec1} can be used to construct the ideal $\fa=(\mathscr{K} \oplus \mathscr{K}^*)^\perp,$ 
and a~decomposition $\fg=\mathscr{K} \oplus \fa\oplus \mathscr{K}^*$, where the generator $x^*$ of $ \mathscr{K}^*$ can be normalized so that ${\mathscr B}_\fg(x,x^*)=1$. 

It has been shown in \cite[Prop. 3.2.2]{BeBou} that $(\fg, {\mathscr B}_\fg)$ is a~double extension of $(\fa, {\mathscr B}_\fa)$, where $ {\mathscr B}_\fa= {\mathscr B}_{\fg}|_{\fa\times \fa}$, by means of a~derivation $\mathscr D\in \fder(\fa)_\od$, $a_0\in \fa_\ev$ and satisfying Eqs. (\ref{2D1}), (\ref{2D2}), (\ref{2D3}).

Let us write $\Delta$ as ($b_0\in \fa_\od$, and $\beta(a), \mu, \nu, \lambda \in \Kee$):
\[
\begin{array}{lcl}
\Delta(x)&=& \lambda x,\\
\Delta(x^*)&=&\nu x^*+b_0+\mu x,\\
\Delta(a)&=& \tilde \Delta(a)+\beta(a)x+\gamma(a) x^*.
\end{array}
\]
Since ${\mathscr B}_\fg(\Delta(x^*),x)={\mathscr B}_\fg(x^*, \Delta(x))$,  then $\nu=\lambda$, and since ${\mathscr B}_\fa(\Delta(a),x)={\mathscr B}_\fa(a, \Delta(x))$,  then $\gamma(a)=0$ for all $a\in \fa$.

In addition,  the conditions 
\[
\begin{array}{l}
{\mathscr B}_\fg(\Delta(a),x^*)={\mathscr B}_\fg ( \tilde \Delta(a)+\beta(a)x ,x^*)=\beta(a) \; \text{ for all $a\in \fa$},\\[2mm]
{\mathscr B}_\fg(a, \Delta(x^*))={\mathscr B}_\fg(a, \lambda x^*+b_0+\mu x)={\mathscr B}_\fa(a,b_0)\; \text{ for all $a$}
\end{array}
\]
imply that 
\begin{equation}\label{eqb1}
\beta(a)={\mathscr B}_\fa(a,b_0) \text{ for all $a\in \fa$}.
\end{equation}
Now, we have 
\[
\begin{array}{lcl}
\Delta([x^*,a]_\fg)&=&\Delta({\mathscr D}(a))=\tilde \Delta ({\mathscr D}(a))+{\mathscr B}_\fa({\mathscr D}(a),b_0)x=\tilde \Delta ({\mathscr D}(a))+{\mathscr B}_\fa(a, {\mathscr D}(b_0))x.
\end{array}
\]
On the other hand, we have
\[
\begin{array}{lcl}
[\Delta(x^*),a]_\fg+[x^*,\Delta(a)]_\fg & = & \lambda {\mathscr D}(a)+[b_0,a]_\fa+{\mathscr B}_\fa({\mathscr D}(b_0),a)x+{\mathscr D}(\tilde \Delta(a)).
\end{array}
\]
Since $\Delta$ is a~derivation, it follows that 
\begin{eqnarray}
\label{eqb2}
[\tilde \Delta, {\mathscr D}] & = &\lambda {\mathscr D}+\ad_{a_0}.
\end{eqnarray}
Following the same computation as in \cite{BaB}, the condition 
\[
\Delta([a,b]_\fg)=[\Delta(a),b]_\fg+[a,\Delta(b)]_\fg
\] 
implies that 
\begin{eqnarray}
\label{eqb3} \tilde \Delta([a,b]_\fa)&=&[\tilde \Delta(a),b]_\fa+[a, \tilde \Delta(b)]_\fa, \\
\label{eqb4} \beta([a,b]_\fa)&=&{\mathscr B}_\fa({\mathscr D}(a),b)+ {\mathscr B}_\fa({\mathscr D}(\tilde \Delta(a)),b)+ {\mathscr B}_\fa({\mathscr D}(a), \tilde \Delta(b)).
\end{eqnarray}
Eq. (\ref{eqb4}) is certainly satisfied, thanks to Eqs. (\ref{eqb1}), (\ref{eqb2}) and the fact that ${\mathscr B}_\fa$ is $\tilde \Delta$-invariant.

Let us consider the condition 
\[
\Delta(s_\fg(f))=[f,\Delta(f)]_\fg \text{ for all $f\in \fg_\od$}.
\]

The LHS ($f=rx+a+tx^*$):
\[
\begin{array}{lcl}
\Delta(s_\fg(f)) & =&\Delta( s_\fa(a) + t^2 a_0+t {\mathscr D}(a))\\[2mm]
&=& \tilde \Delta(s_\fa(a) + t^2 a_0+t {\mathscr D}(a))+{\mathscr B}_\fa(s_\fa(a) + t^2 a_0+t {\mathscr D}(a),b_0)x\\[2mm]
&=& \tilde \Delta(s_\fa(a) + t^2 a_0+t {\mathscr D}(a)).
\end{array}
\]
The RHS ($f=rx+a+tx^*$):
\[
\begin{array}{lcl}
[\Delta(f),f]_\fg& = & [r \lambda x+\tilde \Delta(a)+{\mathscr B}_\fa(a,b_0)x+t(\lambda x^*+b_0+\mu x), rx+a+tx^*]_\fg\\[2mm]
& = & [\tilde \Delta(a)+t(\lambda x^*+b_0), a+tx^*]_\fg\\[2mm]
& = & [\tilde \Delta(a),a]_\fa+t  {\mathscr D}(\tilde \Delta(a))+t \lambda {\mathscr D}(a)+t [b_0,a]_\fa+t^2 {\mathscr D}(b_0).
\end{array}
\]
Taking into account Eq.(\ref{eqb2}),  we see that 
\begin{eqnarray}
\label{eq5} \tilde \Delta (s_\fa(a)) & = & [a,\tilde \Delta (a)]_\fa,\\
\label{eq6} \tilde \Delta(a_0)&=&{\mathscr D}(b_0).
\end{eqnarray}
Eqs. (\ref{eqb3}), (\ref{eq5}) imply that $\tilde \Delta$ is a~derivation of $\fa$; obviously $\tilde \Delta$ is invertible. Define
\[
\omega_\fa(\cdot,\cdot):={\mathscr B}_\fa(\tilde \Delta (\cdot), \cdot).
\]
The form $\omega_\fa$ on $\fa$ is closed,  $\od$-antisymmetric and ortho-orthogonal.\end{proof}
\subsection{Double extensions of Manin triples in chararcteristic $2$}\label{dermanin}
The notion of double extension of Manin triples over a~field  of characteristic zero was introduced in \cite{MR2}. It has been used in  \cite{BBM} to describe the structure of NIS Lie algebras with a closed symplectic form. Later the construction has been superized  in \cite{BaB}. Here this notion is introduced and studied in the case of Lie superalgebras in characteristic $2$, featuring a~condition on both the squaring and the derivation. This unusual feature has already been observed in \cite{BeBou} when dealing with double extensions for Lie superalgebras in characteristic $2$.\\
\begin{mTheorem}[${\mathscr D}_\ev$-extension of Manin triples]\label{MainMn1} Let 
$(\fh, \fg, \fk)$ be a~Manin triple of charac-teristic $2$ with a~NIS ${\mathscr B}_\fh$. Let ${\mathscr D}\in \fder(\fh)_\ev$ satisfy Eq. \eqref{D1}, ${\mathscr D}(\fk_\ev)\subseteq \fk_\ev$, and let $\alpha:\fh_\od\rightarrow \Kee$ be a~quadratic form satisfying Eq. \eqref{D3} and $\alpha(\fk_\od)=0$. Denote by $(\tilde \fh = {\mathscr K} \oplus \fh \oplus {\mathscr K}^*, \tilde {\mathscr B}_{\tilde \fh})$ the double extension of $(\fh, {\mathscr B}_\fh)$
by the one-dimensional Lie algebra ${\mathscr K}$ via ${\mathscr D}$ and $\alpha$.  If ${\mathscr B}_{\tilde \fh}(x^*,x^*)=0$, then the triple $(\tilde \fh, \tilde \fg, \tilde \fk)$ is also a~Manin triple with a~NIS  $\tilde {\mathscr B}_{\tilde \fh}$,  where $\tilde \fg= {\mathscr K} \oplus \fg$ and $\tilde \fk= {\mathscr K}^* \ltimes \fk$.
 
\end{mTheorem}

We call the Manin triple $(\tilde \fh, \tilde \fg, \tilde \fk)$ with NIS $\tilde {\mathscr B}_{\tilde \fh}$ {\it the double extension of the Manin triple $(\fh, \fg, \fk)$ with NIS ${\mathscr B}_{\fh}$ by the one-dimensional Lie algebra ${\mathscr K}$}.

\begin{proof} Let us first show that if $\alpha(\fk_\od)=0$,  then ${\mathscr D}(\fk_\od)\subseteq \fk_\od$. Indeed, let ${\mathscr D}(k)\in {\mathscr D}(\fk_\od)$ and write ${\mathscr D}(k)=\mathrm{proj}_{\fg}({\mathscr D}(k))+\mathrm{proj}_{\fk}({\mathscr D}(k))$. For all $l\in \fk_\od$, we have \footnote{Recall that ${\mathscr B}_\alpha$ is the polar form of the quadratic form $\alpha$.}
\[
0={\mathscr B}_\alpha(k,l)={\mathscr B}_\fh({\mathscr D}(k),l)={\mathscr B}_\fh(\mathrm{proj}_{\fg}({\mathscr D}(k)),l).
\]
On the other hand, since ${\mathscr B}_\fh$ is even and ${\mathscr B}_\fh(\fg,\fg)=0$,  it follows that 
\[
{\mathscr B}_\fh(\mathrm{proj}_{\fg}({\mathscr D}(k)),l)=0 \; \text{ for all $l\in \fh$}.
\]
It follows that $\mathrm{proj}_{\fg}({\mathscr D}(k))=0$ since ${\mathscr B}_\fh$ is non-degenerate.

Let us show that $\tilde \fg=  \fg \oplus {\mathscr K}$ is Lie sub-superalgebra of $\tilde \fh$. Indeed, for every $\lambda x+f, \mu x+g\in \tilde \fg$, where $\lambda,\mu \in \Kee$ and $f,g\in \fg$, we have 
\[
[\lambda x+f, \lambda x+g]_{\tilde \fh}=[f,g]_{\tilde \fh}=[f,g]_\fh+{\mathscr B}_\fh(f,g)x \in  \tilde \fg,
\]
since $\fg$ is a~Lie sub-superalgebra of $\fh$. Now the squaring; for all $f\in \tilde \fg_\od$, we have
\[
s_{\tilde \fh}(f)=s_\fh(f)+\alpha(f)x  \in  \tilde \fg_\ev,
\]
since $\fg$ is a~Lie sub-superalgebra of $\fh$. 

Let us show that $\tilde \fk= {\mathscr K}^* \ltimes \fk$ is a Lie sub-superalgebra of $\fh$.  For every ${\lambda x^*+f, \mu x^*+g\in \tilde \fk}$, where $\lambda,\mu \in \Kee$ and $f,g\in \fk$, we have 
\[
\begin{array}{lcl}
[\lambda x^*+f,\mu x^*+g]_{\tilde \fh} & = & [f,g]_\fh+{\mathscr B}_{\fh}({\mathscr D}(f),g)x+\mu {\mathscr D}(f)+ \lambda {\mathscr D}(g)  \in  \tilde \fk\\[2mm]
& = & [f,g]_\fh+\mu {\mathscr D}(f)+ \lambda {\mathscr D}(g)  \in  \tilde \fk.
\end{array}
\]
since $\fk$ is a~Lie subsuperalgebra of $\fh$  and ${\mathscr B}_\fh({\mathscr D}(f),g)=0$ since ${\mathscr D}(\fk)\subseteq \fk$. Now the squaring; for all $f\in \tilde \fk_\od$, we have
\[
s_{\tilde \fh}(f)=s_\fh(f)+\alpha(f)x =s_\fh(f) \in  \tilde \fk_\ev,
\]
since $\fk$ is a~Lie subsuperalgebra of $\fh$ and $\alpha(\fk_\od)=0$.

The condition ${\mathscr B}_{\tilde \fh}(\tilde \fg, \tilde \fg)={\mathscr B}_{\tilde \fh}(\tilde \fk, \tilde \fk)=0$ follows from the explicit construction of ${\mathscr B}_{\tilde \fh}$ as in Theorem \ref{MainTh} and ${\mathscr B}_{\tilde \fh}(x^*,x^*)=0$.
\end{proof}

Now, we will give the converse of Theorem \ref{MainMn1}.\\

\begin{mTheorem}[Converse of Theorem \ref{MainMn1}]\label{MainMnR1}
Let $(\fh, \fg,  \fk)$ be a~Manin triple such that $\dim(\fh)>1$.  If  $(\fz_s(\fh)\cap \fg)_\ev \not =0$ or $(\fz_s(\fh)\cap \fk)_\ev \not =0$, then the Manin triple $(\fh, \fg, \fk)$ is a~double extension of a~Manin triple $(\fc, \fa, \fb)$ as in Theorem $\ref{MainMn1}$.


\end{mTheorem}
\begin{proof} 

Let us assume that $(\fz_s(\fh)\cap \fg)_\ev \not =0$. The case $(\fz_s(\fh)\cap \fk)_\od \not =0$ is absolutely the same. Let us choose a non-zero $x\in (\fz_s(\fh)\cap \fg)_\ev$ and write ${\mathscr K}=\Span\{x\}$. By construction, ${\mathscr K}$ and ${\mathscr K}^\perp$ are both ideals. Since ${\mathscr B}_\fh(\fg,\fg) = \{0\}$,  we obtain $\fg \subseteq {\mathscr K}^\perp$.  Since ${\mathscr B}_\fh$ is non-degenerate and ${\mathscr B}_\fh(\fg,\fg) = \{0\}$,   then there exists $x^*  \in \fk_\ev$ such that
${\mathscr B}_\fh(x^*,x) \not=0$. We may assume that ${\mathscr B}_\fh(x^*,x)=1$. Let us write ${\mathscr K}^*=\Span\{x^*\}$. Since $\mathscr{K}\subset \mathscr{K}^\perp$, it follows that $\dim(\mathscr{K}^\perp)=\dim(\fh)-1$. We can then write $\fh={\mathscr K}^\perp \oplus {\mathscr K}^*$ since $x^*\not \in {\mathscr K}^\perp$. There exists a vector subsuperspace $\fc:=({\mathscr K}\oplus{\mathscr K}^*)^\perp \subseteq \fh$ such that ${\mathscr K}^\perp=  {\mathscr K}\oplus \fc$ and ${\mathscr B}_\fh(\fc, x^*)=\{0\}$. We easily see that ${\mathscr B}_{\fc}:={\mathscr B}_{\fh}|_{\fc\times \fc}$ is non-degenerate. Since $\fg \subseteq {\mathscr K}^\perp$ we have $\fg={\mathscr K} \oplus \fa$, where $\fa=\fg \cap \fc$. 

On the other hand, as $\fh=\fg\oplus \fk$ and $\fg \subseteq {\mathscr K}^\perp$ then ${\mathscr K}^\perp=\fg\oplus \fb$,  where $\fb=\fk\cap {\mathscr K}^\perp$.  Since ${\mathscr B}_\fh(\fk,\fk)=\{0\}$, ${\mathfrak c}=({\mathscr K}\oplus {\mathscr K}^*)^\perp$ and ${\mathfrak b}={\mathfrak k} \cap {\mathscr K}^\perp$, we get $\fb \subseteq \fc$. Then, $\fh={\mathscr K}\oplus \fc \oplus {\mathscr K}^*$ and $\fc=\fa \oplus \fb$ by the dimension considerations.  Recall that $\fc\subseteq {\mathscr K}^\perp$. We get
\[
[f,g]_\fh=\beta(f,g)+\phi(f,g)x \text{ for all $f,g\in \fc$},
\]
where $\beta(f,g)\in \fc$ and $\phi(f,g)\in  \Kee$. Moreover, 
\[
[x^*,f]_\fh={\mathscr D}(f)+\psi(f)x \text{ for all $f\in \fc$},
\]
with ${\mathscr D}(f)\in \fc$ and $\psi(f)\in \Kee$. 

Moreover, the squaring reads as follows
\[
s_{\fh}(f)=s_{\fc}(f)+\alpha(f) x \text{ for all $f\in \fc_\od$,}
\]
where $s_\fc(f)\in \fc_\ev$.

Now, let $f\in \fb_\od \subset \fk_\od$. Since $s_\fh(f)\in \fk_\ev$,  then ${\mathscr B}_\fh(s_\fh(f),x^*)=0$ as ${\mathscr B}_\fh(\fk, \fk)=0$. It follows that 
\[
0={\mathscr B}_\fh(s_\fh(f),x^*)={\mathscr B}_\fh (s_\fc(f) + \alpha(f)x,x^*)=\alpha(f).
\]
Therefore, $\alpha(\fb_\od)=\{0\}$.

Following the same proof as in \cite[Prop. 3.1.3]{BeBou}, we can show that $(\fc, \beta, s_\fc, {\mathscr B}_\fc)$ is a~double extension of $(\fh, [\cdot,\cdot]_\fh, s_\fh, {\mathscr B}_\fh)$ by means of the derivation ${\mathscr D}$ and the quadratic form $\alpha$. The proof also implies that $\psi=0$ and that 
\[
\phi(f,g)={\mathscr B}_\fc({\mathscr D}(f),g)\text{ for all $f,g\in \fc$}.
\]

Let us show that ${\mathscr D}(\fb_\ev)\subseteq \fb_\ev$. Since $\fb=\fk \cap {\mathscr K}^\perp$, but ${\mathscr K}^\perp$ is an ideal in $\fh$ and $\fk$ is a~sub-superalgebra of $\fh$,  it follows that for all $f\in \fb_\ev$, we have
\[
{\mathscr D}(f) =[x^*,f] \in(\fk \cap {\mathscr K}^\perp)_\ev= \fb_\ev.
\]
It remains to show that ${\mathscr B}_\fc(\fa,\fa)={\mathscr B}_\fc(\fb,\fb)=0$,  but this is obvious since $\fa\subset \fg$ and $\fb\subset \fk$.\end{proof}

Hereafter,  we will carry out the study in the case of $\mathscr D_\od$-extensions. \\

\begin{mTheorem}[${\mathscr D}_\od$-extension for Manin triples]\label{MainMn2} Let $(\fh, \fg, \fk)$ be a~Manin triple in charac-teristic $2$ with a~NIS ${\mathscr B}_\fh$. Let ${\mathscr D}\in \fder(\fh)_\od$ satisfy Eqns. \eqref{2D1}, \eqref{2D2}, \eqref{2D3} and ${\mathscr D}(\fk) \subseteq \fk$. Let $(\tilde \fh = {\mathscr K} \oplus \fh \oplus {\mathscr K}^*, \tilde {\mathscr B}_{\tilde \fh})$ be the double extension of $(\fh, {\mathscr B})$
by the one-dimensional Lie superalgebra ${\mathscr K}$ (by means of ${\mathscr D}$). If $a_\ev \in \fk_\ev$, then $(\tilde \fh, \tilde \fg, \tilde \fk)$ is a~Manin triple with a~NIS $\tilde {\mathscr B}_{\tilde \fh}$,  where $\tilde \fg= {\mathscr K} \oplus \fg$ and $\tilde \fk= {\mathscr K}^* \ltimes \fk$. 

\end{mTheorem}
\begin{proof}

Similar to that of Theorem \ref{MainMn1}.
\end{proof}

Now, let us give the converse of Theorem \ref{MainMn2}. \\

\begin{mTheorem}[Converse of Theorem \ref{MainMn2}]\label{MainMnR2}
Let $(\fh, \fg,  \fk)$ be a~Manin triple such that $\dim(\fh)>1$.  If  $(\fz(\fh)\cap {\mathscr C}(\fh, {\mathscr B}_\fh)\cap \fg)_\od \not =0$ or $(\fz(\fh)\cap {\mathscr C}(\fh, {\mathscr B}_\fh)\cap \fk)_\od \not =0$ then the Manin triple $(\fh, \fg, \fk)$ is a~double extension of a~Manin triple $(\fc, \fa, \fb)$ as in Theorem $\ref{MainMn2}$.


\end{mTheorem}
\begin{proof}

Similar to that of Theorem \ref{MainMnR1}.
\end{proof}

\section{Conclusion and outlook}

(1) Manin triples for Lie superalgebras in characteristic $2$ are studied and two  examples are given in \S \ref{examples}. The notion of Lie bi-superalgebra has also been studied and a link has been established between the two.  Because of the squaring, the conditions for the existence of a Lie bi-superalgebra structure are more complex than for characteristic not 2.

(2) The construction has been generalized to any pair of Lie superalgebras $\fg$ and $\fk$ satisfying certain conditions.  Such a pair satisfying these conditions is called  a~\textit{matched pair}.

(3) It has been showed that if there is a~classical $r$-matrix on the~Lie algebra $F(\fg)$,  where $\fg$ is a Lie superalgebra and $F$ is the functor of forgetting the super structure,  then a~squaring on $\fg^*$ can be defined.  It has been showed  that this squaring satisfies the Jacobi identity if and only if a~certain condition \eqref{JIsqdual} is satisfied. Such classical $r$-matrices are called {\it admissible}.  Surprisingly, having a classical admissible $r$-matrix guarantees the existence of Lie bi-superalgebra on $\fg$ in characteristic 2. 

{\bf Open problem:} Following \cite{F1, F2}, which Lie superalgebras in characteristic 2 possess at least one admissible non-trivial classical matrix $r$? At least describe such $r$-matrices in the case of known Lie superalgebras in characteristic 2 (see \cite{BGL}, \cite{BGLLS}).

(4) In particular,  admissible classical $r$-matrices on Lie superalgebras in characteristic 2 with an~even NIS have been studied at the end of \S\ref{rmatManin}.  It has been showed that every such Lie superalgebras endowed additionally with a closed $\od$-antisymmetric ortho-orthogonal form posses an admissible classical $r$-matrix.  More precisely,  these $r$-matrices are obtained by the inverse of the derivations that links NIS and the ortho-orthogonal form.  We have developed a process of double extensions to describe these $r$-matrices,  see \S\ref{DEeven}.  

{\bf Open problems:} (i) Following \cite{BBM},  is it possible to obtain every Lie superalgebra $\fg$ in characteristic 2 with an even NIS and a closed $\od$-antisymmetric ortho-orthogonal form as a double extension of a smaller Lie superalgebra $\fa$ (with the same property) such that~$\dim(\fg)=\dim(\fa)+2$? (ii) Study this problem in the case where NIS is odd.  We believe that this can be done by studying Manin triples of the form $(\fg\oplus \Pi(\fg^*),  \fg, \fg^*)$, where $\Pi$ is the functor of changing the parity.

(5) It has been proved in \cite{BaB} that every even NIS-Lie superalgebra equipped with a closed antisymmetric ortho-symplectic form over an algebraically closed field of characteristic zero is actually a~Manin triple. The proof makes use of the weight space decomposition of a~ Lie superalgebra that is only valid over a~field of characteristic zero. We do not know whether the result of \cite{BaB} is true or not true in our context, perhaps ideas from \cite{BLLS} could help.

(6) The importance of invertible derivations on Lie superalgebras in characteristic 2 (not necessarily NIS) is crucial to construct admissible Yang-Baxter operators,  as presented in Example 3 of \S \ref{examples} and Example 4 in \S \ref{DEeven}. It is a well-known theorem by Jacobson \cite{J} that every Lie algebra in characteristic 0 with an invertible derivation has to be nilpotent.

{\bf Open problem:} Is Jacobson's Theorem still valid for Lie (super)algebras in characteristic 2? If the theorem doesn't hold describe the structure of such Lie (super)algebras.


(7) It has been shown that every Lie algebra with an invertible derivation gives rise to a~left-symmetric structure on it.  Further, it has been shown in Theorem \ref{WNIS} that every ${\mathbb Z}/2$-graded Lie algebra in characteristic $2$ with an invertible derivation can be turned into a Lie superalgebra provided that the left-symmetric structure is left-alternative. 

{\bf Open problems:} (i) It has been proved that every left-alternative structure is left-symmetric, see Prop.  \ref{leftalt}.  What about the converse? (ii) If the converse doe not hold,  describe invertible derivations on ${\mathbb Z}/2$-graded Lie algebras in characteristic 2 for which the left-symmetric structure is left-alternative.  (iii) Compare the method 2 in \cite{BLLSq} with our method in constructing Lie superalgebras from ${\mathbb Z}/2$-graded Lie algebras in characteristic 2.\\

{\bf Acknowledgment:} We are very thankful to the referee for providing comments that helped us improve the paper's presentation.



\begin{thebibliography}{999}
\bibitem[ABBQ]{ABBQ} Albuquerque H., Barreiro E. and Benayadi S., Odd quadratic Lie superalgebras. J. of Geometry and Physics, {\bf 60} (2010), 230-250.

\bibitem[ABB]{ABB} Albuquerque H., Barreiro E. and Benayadi S., Quadratic Lie superalgebras with a~reductive even part. J. Pure Appl. Algebra, {\bf 213} (2009), 724--731.

\bibitem[BBB]{BBB}  Bajo I., Benayadi S. and Bordemann M., Generalized double extension and descriptions of quadratic Lie superalgebras, \texttt{arXiv:0712.0228}.

\bibitem[BBM]{BBM} Bajo I., Benayadi S. and Medina A., Symplectic structures on quadratic Lie algebras. Journal of Algebra, {\bf 316} (2007), 174--188.

\bibitem[BaB]{BaB} Barreiro E. and Benayadi S., Quadratic symplectic Lie superalgebras and Lie bi-superalgebras. Journal of Algebra, {\bf 321} (2009), 582--608.

\bibitem[BD1]{BD1} Belavin A. A. and Drinfeld V. G., Solutions of the classical Yang-Baxter equation and simple Lie algebras. Funct. Anal. Appl. 16 (1982), 159--180.

\bibitem[BD2]{BD2} Belavin A. A. and Drinfeld V. G., Triangle equation and simple Lie algebras. Soviet Sci. Reviews C4 (1984), 93--165.


\bibitem[BeB]{BeB} Benamor H. and Benayadi S., Double extension of quadratic Lie superalgebras. Comm. Algebra, {\bf 27}, No. 1 (1999), 67--88.

\bibitem[B]{B} Benayadi S., Quadratic Lie superalgebras with completely reductive action of the even part on the odd part. J. of Algebra, {\bf 223} (2000), 344-366.

\bibitem[BeBou]{BeBou} Benayadi S. and Bouarroudj S., Double extensions of Lie superalgebras in characteristic $2$ with non-degenerate invariant supersymmetric bilinear forms. J. of Algebra, {\bf 510} (2018), 141--179; arXiv:1707.00970.

\bibitem[BBH]{BBH} Benayadi S., Bouarroudj S., Hajli M., Double extensions of restricted Lie (super)algebras. Arnold. Math. J.  {\bf 6} (2020), 231 -- 269; \texttt{arXiv:1810.03086}

\bibitem[Bor]{Bor} Bordemann M., Nondegenerate invariant bilinear forms on nonassociative algebras. Acta Math. Univ. Comenianae, Vol. LXVI, 2 (1997), 151--201.



\bibitem[BGL]{BGL} Bouarroudj S., Grozman P., Leites D., Classification of finite-dimensional modular Lie superalgebras with indecomposable Cartan matrix. Symmetry, Integrability and Geometry: Methods and Applications (SIGMA), 5 (2009), 060, 63 pages; arXiv:math.RT/0710.5149

\bibitem[BGLLS]{BGLLS}
Bouarroudj S., Grozman P., Lebedev A., Leites D. and Shchepochkina I., Simple vectorial Lie algebras in characteristic
$2$ and their superizations. Symmetry, Integrability and Geometry: Methods and Applications (SIGMA) {\bf 16} (2020), 089, 101 pages; \texttt{arXiv:1510.07255}





\bibitem[BLLS]{BLLS} Bouarroudj S., Leites D., Lozhechnyk O. and Shang J. The roots of exceptional modular Lie superalgebras with Cartan
matrix. Arnold Math. J.  {\bf 6} (2020), 63 --118; \texttt{arXiv:1904.09578}

\bibitem[BGL2]{BGL2} Bouarroudj S., Grozman P., and Leites D., Deforms of symmetric simple modular Lie algebras and Lie superalgebras; \texttt{arXiv:0807.3054}

\bibitem[BLLSq]{BLLSq} Bouarroudj S., Lebedev A., Leites D. and Shchepochkina I., Classifications of simple Lie superalgebras in characteristic $2$. Int.  Math.  Res.  Not. (2021); https://doi.org/10.1093/imrn/rnab265.

\bibitem[BLS]{BLS} Bouarroudj S., Leites D. and Shang J., Computer-aided study of double extensions of restricted Lie superalgebras preserving the non-degenerate closed 2-forms in characteristic $2$. Experimental Mathematics (2019); https://doi.org/10.1080/10586458.2019.1683102.

\bibitem[CP]{CP} Chari, V. and Pressley, A., A guide to quantum groups. Cambridge University Press, Cambridge, 1994.

\bibitem[D]{D} Drinfeld V. G., Quantum groups. J. Sov. Math. {\bf 41} (1988), 898--915.

\bibitem[FS]{FS} 
Favre G. and Santharoubane L. J., Symmetric, invariant, non-degenerate bilinear form on a~Lie algebra. J. of Algebra, {\bf 105} (1987), 451--464. 

\bibitem[F1]{F1} Feldvoss J., Existence of triangular Lie bialgebra structures. J. Pure Appl. Algebra {\bf 134} (1999), 1--14.

\bibitem[F2]{F2} Feldvoss J., Existence of triangular Lie bialgebra structures II. J. Pure Appl. Algebra {\bf 198} (2005), 151--163.

\bibitem[GZB]{GZB} Gould M. D., Zhang R. B.,  Bracken A. J., Lie bi-superalgebras and graded classical Yang–Baxter equation. Rev. Math. Phys. {\bf 3} (1991) 223--240.

\bibitem[J]{J} Jacobson N.,  A note on automorphisms and derivations of Lie algebras,  Proc. Amer. Math. Soc. 6 (1955) 33--39.

\bibitem[Ka1]{Ka1} Karaali G., Constructing $r$-matrices on simple Lie superalgebras. J. of Algebra {\bf 282} (2004), 8-102.

\bibitem[Ka2]{Ka2} Karaali G.,  A new Lie bialgebra structure on $\mathrm{sl}(2,1)$.  Representations of algebraic groups, quantum groups,  and Lie algebras, 101--122,  Contemp.  Math.,  413,  Amer. Math. Soc.,  Providence,  RI,  2006.

\bibitem[LeD]{LeD} Lebedev A., Analogs of the orthogonal, Hamiltonian, Poisson, and contact Lie superalgebras in
characteristic $2$. J. Nonlinear Math. Phys. 17 (2010), suppl. 1, 217--251.



\bibitem[LSoS]{LSoS}
Leites D. (ed.) \textit{Seminar on supersymmetry v. $1$. Algebra and
Calculus: Main chapters}, (J.~Bernstein, D.~Leites, V.~Molotkov,
V.~Shander), MCCME, Moscow,  (2012),  410 pp (in Russian; a~version in
English is in preparation but available for perusal)

\bibitem[Lyb]{Lyb}
Leites D. Corrections and additions to the paper ``D.~Leites, V.~Serganova. Solutions of the classical Yang--Baxter equation
for simple Lie superalgebras. 
Theor. and Mathem. Phys. \textbf{58} (1984), no.~1, 16--24. MR
85m:82020.'' 

\bibitem[LSh]{LSh} Leites D. and Shapovalov A., Manin-Olshansky triples and Lie superalgebras. J. Nonlinear Math.
Phys., v. 7, 2000, no. 2, 120--125; arXiv:math.QA/0004186

\bibitem[LS]{LS} Leites D. and Serganova V., Solutions of the classical Yang-Baxter equation for simple superalgebras. Theoretical and Mathematical Physics {\bf 58} (1984), 16--24.

\bibitem[Maj]{Maj} Majid S., Matched pairs of Lie groups associated to solutions of the Yang-Baxter equations. Pacific J. Math. 141(2): (1990) 311--332.

\bibitem[MR1]{MR1} Medina A. and Revoy Ph., Alg\`ebres de Lie et produit scalaire invariant, Ann. Sci. \'Ecole Norm. Sup. (4) 18 (1985) 553--561.

\bibitem[MR2]{MR2} Medina A. and Revoy Ph., Le notion de double extension et les groupes de Lie-Poisson. Semin. Gaston Darboux Geom. Topologie Differ. 1987--1988 (1988) 141--171.

\bibitem[M]{M} Michaelis W., A class of infinite-dimensional Lie bialgebras containing the Virasoro algebra, Adv. in Mathematics {\bf 107} (1994), 365--392.

\bibitem[M2]{M2} Michaelis W.,  Lie coalgebras.  Adv. Math. {\bf 38} (1980),  no. 1,  1--54. 

\bibitem[N]{N} Nijenhuis A., Sur une classe de propri\'et\'es communes \`a quelques types diff\'erents d'alg\`ebres.  L'enseignement Mathem. {\bf XIV}, Fasc. 3-4,  (1968),  225--275. 

\bibitem[O]{O} Olshansky G., Quantized universal enveloping superalgebra of type Q and a~superextension of the Hecke algebra. Lett. Math. Phys., (1992), V.24, N 2, 93--102.





\bibitem[Sm]{Sm} de Smedt V., Existence of a~Lie bialgebra structure on every Lie algebra. Lett. Math. Phys. {\bf 31} (1994), 225--231.

\bibitem[V]{V} Vinberg E. B., Theory of convex homogeneous cones. Trydu Moscow Mat. Obshch. 12 (1963) 303--358; Transl. Moscow Math. Soc. 12 (1963) 340--403.

\end{thebibliography}
\end{document}